\documentclass[12pt]{article}
\usepackage[margin=1in]{geometry}

\usepackage{mathptmx}

\usepackage{amsthm,amsmath,amsfonts,amssymb}
\usepackage[numbers,sort&compress]{natbib}
\usepackage{microtype}
\usepackage{graphicx}
\usepackage{float}
\usepackage{needspace}
\usepackage{enumitem}
\usepackage[colorlinks,citecolor=blue,urlcolor=blue,linkcolor=blue]{hyperref}
\hypersetup{
  pdftitle={De Finetti + Sanov = Bayes: Exchangeable Prediction under Moment Constraints},
  pdfauthor={Nicholas G. Polson and Daniel Zantedeschi},
  pdfsubject={Exchangeable prediction under moment constraints}
}

\theoremstyle{plain}
\newtheorem{theorem}{Theorem}
\newtheorem{proposition}{Proposition}
\newtheorem{corollary}{Corollary}
\newtheorem{lemma}{Lemma}
\theoremstyle{definition}

\newtheorem{remark}{Remark}
\newtheorem{example}{Example}

\newcommand{\X}{\mathcal{X}}
\newcommand{\PN}{\widehat{P}_N}
\newcommand{\KL}[2]{D(#1\,\|\,#2)}
\newcommand{\Ber}{\mathrm{Ber}}
\newcommand{\Ex}{\mathbb{E}}
\renewcommand{\Pr}{\mathbb{P}}
\newcommand{\Qbb}{\mathbb{Q}}
\newcommand{\Real}{\mathbb{R}}

\newcommand{\Pstar}{P^{\star}}
\newcommand{\Fish}{\mathcal{I}}
\newcommand{\TN}{\mathcal{T}_N}
\newcommand{\simplex}{\Delta}
\newcommand{\TV}[1]{\lVert #1\rVert_{\mathrm{TV}}}
\newcommand{\PiE}{\widetilde{\Pi}_E}
\DeclareMathOperator*{\argmin}{arg\,min}

\let\origparagraph\paragraph
\renewcommand{\paragraph}[1]{\origparagraph{#1.}}

\let\origappendix\appendix
\renewcommand{\appendix}{\clearpage\origappendix}

\title{De Finetti + Sanov = Bayes:\\
Exchangeable Prediction under Moment Constraints}

\author{Nicholas G. Polson\thanks{Booth School of Business, University of
Chicago, \texttt{ngp@chicagobooth.edu}.}
\and Daniel Zantedeschi\thanks{Muma College of Business, University of South
Florida, \texttt{danielz@usf.edu}.}}

\date{\today}

\begin{document}

\maketitle

\begin{abstract}
We study exchangeable prediction when empirical-moment constraints define each active
finite horizon $N$. The relevant law is the de Finetti mixture conditioned on
$E_N=\{\widehat{P}_N\in E_{\varepsilon_N}\}$. By permutation invariance, the target
may be any fixed block of $m$ coordinates within the active horizon, regardless of
whether those coordinates are labeled past, held out, or future relative to any
finite cut. Conditionally on the directing measure $\mu$, the Gibbs-conditioning
principle sends the law of such a
block to the $m$-fold product $(\Pstar_\mu)^{\otimes m}$ of the $I$-projection
$\Pstar_\mu=\argmin_{Q\in E}\KL{Q}{\mu}$. On a finite alphabet we give an
elementary master inequality for general polyhedral moment windows. After mixing
over the constraint posterior $\Pi_{N,E}$, and under weak convergence plus
posterior-averaged component control, the finite-dimensional marginals converge to
a consistent exchangeable law whose random directing measure is the $I$-projection
$\Pstar_\mu$, with $\mu$ drawn from the weak limit $\Pi_E$. Sequential prediction
under this limiting law is therefore Bayesian prediction from a mixture of
componentwise $I$-projections. The limiting behavior depends on reachability. For a
reachable constraint, the projection is asymptotically the identity on the selected
subfamily. Under additional regularity, an unreachable constraint makes the
constraint posterior concentrate on the rate-minimizing subfamily, while the
projections remain nontrivial.
In our examples, at least one operation is asymptotically inactive, though both
enter the finite-horizon construction. The
master bound also reads as an equivalence of ensembles. We reserve ``maximum entropy'' for a uniform or flat baseline and use
``minimum relative entropy'' or ``$I$-projection'' in general.
\end{abstract}

\noindent\textbf{Keywords:} Exchangeability; Gibbs conditioning; I-projection; Equivalence of ensembles; Predictive inference.

\thispagestyle{empty}
\pagenumbering{roman}
\newpage
\pagenumbering{arabic}

\section{Introduction}
\label{sec:intro}

De Finetti's theorem and the Gibbs conditioning principle are usually told as
separate stories. The first says that an exchangeable sequence is a mixture of
i.i.d.\ laws indexed by a random directing measure. This is the Bayesian
representation: the mixing measure is the prior, and its conditional law given the
data is the posterior \citep{definetti,heathsudderth,fortinipetrone}.
The second, the conditional form of Sanov's theorem, says that once one conditions
on the empirical measure of an i.i.d.\ sample lying in a constraint set, the law
of a fixed block of coordinates converges to a product of exponential tilts, the
$I$-projection of the sampling law onto the set \citep{csiszar84,vancampenhout,leonardnajim}.
This paper treats the empirical-moment constraint as a primitive feature of the
predictive model and shows how the two theorems combine: Gibbs conditioning
transforms the directing measure, and de Finetti supplies Bayesian prediction
under the transformed law.

\subsection{The constrained-horizon experiment}
\label{sec:experiment}

Fix a finite alphabet $\X=\{1,\dots,k\}$.
We do not begin with a completed
unconstrained sample and subsequently observe a rare event about it. Instead, for
each active prediction horizon $N$, we define the model separately by conditioning
the exchangeable law on an empirical-moment restriction. Start from the de Finetti joint law
$\Pr_N(d\mu,dx_{1:N})=\Pi(d\mu)\,\mu^{\otimes N}(dx_{1:N})$, with directing random
measure $\mu\sim\Pi$, and impose the primitive constraint
\begin{equation}
\label{eq:eventN}
E_N=\{\widehat{P}_N\in E_{\varepsilon_N}\},\qquad \widehat{P}_N=\frac1N\sum_{i=1}^N\delta_{X_i},
\end{equation}
where $E_{\varepsilon_N}$ is a shrinking window around a moment-constraint set $E$
(Section~\ref{sec:windows}). The constrained-horizon joint law is
\begin{equation}
\label{eq:QN}
\Qbb_N(d\mu,dx_{1:N})=\Pi_{N,E}(d\mu)\,\Qbb_{N,\mu}(dx_{1:N}),
\end{equation}
where the \emph{constraint posterior} (or constraint-updated mixing law) is
\begin{equation}
\label{eq:constraintpost}
\Pi_{N,E}(d\mu)=\frac{\Pr_\mu(E_N)\,\Pi(d\mu)}{\int\Pr_\nu(E_N)\,\Pi(d\nu)},
\end{equation}
while $\Qbb_{N,\mu}(\cdot)=\Pr_\mu(\cdot\mid E_N)$ is the component law after
conditioning. This factorization separates two changes. Within a component, the
event alters the joint law of $X_{1:N}$; across components, it favors those values
of $\mu$ for which $E_N$ is more likely. The event itself concerns the whole
horizon and is permutation invariant. It can encode aggregate side information,
a conservation or calibration requirement, or a published sample moment.
Prediction concerns only a few coordinates. Whether those coordinates are called
past, held out, or future is an interpretation imposed outside the exchangeable
law and does not change their marginal conditional distribution.

\paragraph*{Projective consistency}
The laws $\Qbb_N$ need not form a projective family. Restricting
$\Qbb_{N+1}$ to $\X^N$ generally does not recover $\Qbb_N$, since $E_{N+1}$
and $E_N$ condition different horizons. We therefore do not identify the
finite-$N$ laws with the marginals of a single process. The infinite predictive
object is instead obtained from a consistent \emph{local limit} of their
finite-dimensional marginals, much as an infinite-volume Gibbs state is obtained
from finite constrained systems. Section~\ref{sec:limit} carries this out while
keeping the constraint active at every horizon.

\subsection{Local prediction under a global constraint}
\label{sec:local}

Let $J_N,I_N\subseteq\{1,\dots,N\}$ be disjoint index sets with $|J_N|=t$ and
$|I_N|=m$ fixed, $J_N$ an observed set and $I_N$ a prediction set; their locations
may vary with $N$. We study $\Qbb_N(X_{I_N}\in\cdot)$ and, when
observations $X_{J_N}=x_J$ are supplied,
$\Qbb_N(X_{I_N}\in\cdot\mid X_{J_N}=x_J)$. Since $\mu^{\otimes N}$ and $E_N$ are
permutation invariant, only the cardinalities and observed values matter, not the
coordinate labels:
\begin{equation}
\label{eq:locfree}
\mathcal{L}_{\Qbb_{N,\mu}}(X_{I_N})=\mathcal{L}_{\Qbb_{N,\mu}}(X_{1:m}).
\end{equation}
A target labeled future remains subject to the tilt whenever it lies inside the
horizon on which $E_N$ is imposed. Thus temporal location is irrelevant; what
matters is membership in the constrained horizon.

\subsection{A location-free Bernoulli example}
\label{sec:bernoulli-intro}

\begin{example}[Location-free Bernoulli prediction]
\label{ex:bernoulli}
Let the baseline component be $\mu=\Ber(0.2)$ and, for each horizon $N$, impose
the primitive condition $\bar X_N\approx0.8$. The $I$-projection of $\Ber(0.2)$
onto the mean-$0.8$ constraint is $\Pstar_\mu=\Ber(0.8)$, so for any sequence of
indices $i_N\le N$,
\[
\mathcal{L}_{\Qbb_{N,\mu}}(X_{i_N})\ \Longrightarrow\ \Ber(0.8).
\]
The index $i_N$ may be an unrevealed early coordinate, a held-out coordinate, or a
coordinate labeled future relative to an arbitrary cut $t<i_N$; since it belongs
to the active constrained horizon, its local law is the same. More generally, for
fixed $t$ and a history $x_{1:t}$ with positive limiting probability,
$\Qbb_{N,\mu}(X_{t+1}=1\mid X_{1:t}=x_{1:t})\to0.8$.
\end{example}

\subsection{Scope}
\label{sec:scope}

De Finetti alone supplies the Bayesian representation, and Sanov is not needed to
derive Bayes' rule. The role of Gibbs conditioning is to identify the component
law generated by the primitive empirical constraint. The constrained finite-horizon
ensembles converge locally to an exchangeable law whose random directing measure
is the componentwise $I$-projection $\Pstar_\mu$; de Finetti then supplies ordinary
Bayesian prediction under that limiting constrained law. The title's ``$=$''
refers to this transformation of the directing law, not to a literal addition of
theorems.

\subsection{Prior art and contribution}
\label{sec:contribution}

Conditional exponential tilting and equivalence of ensembles are classical, and one
identity motivates the whole construction. For $X_i$ i.i.d.\ $\mu$ conditioned on
an exact empirical average of $g$, \citet{vancampenhout} show, under their regularity
conditions and along horizons for which the conditioning event is attainable, that
a single coordinate has limit law
\begin{equation}
\label{eq:vcc}
\begin{aligned}
&\Pr_\mu\Big(X_1=x\ \Big|\ \tfrac1N\textstyle\sum_{i=1}^N g(X_i)=\alpha\Big)\\
&\qquad\longrightarrow\ \Pstar_\mu(x)
=\frac{\mu(x)\,e^{\lambda^\top g(x)}}{\sum_{y\in\X}\mu(y)\,e^{\lambda^\top g(y)}},
\end{aligned}
\end{equation}
with $\lambda$ fixed by $\langle g,\Pstar_\mu\rangle=\alpha$: the conditioned
coordinate is drawn not from the baseline $\mu$ but from its exponential tilt, which
is exactly the $I$-projection $\argmin_{Q\in E}\KL{Q}{\mu}$ of \eqref{eq:iproj}.
The conditioned coordinates are not independent at finite $N$, since
conditioning on the empirical measure couples them; independence is recovered only
in the limit. \citet{csiszar84} supplies the corresponding block statement,
convergence to the product $(\Pstar_\mu)^{\otimes m}$, in relative-entropy geometry.
The identity is used in statistics: \citet{sterncover} recover the distribution of
lottery-ticket combinations from observed marginals against a uniform baseline, where
maximum entropy is the correct description, and projection onto marginal constraints
is iterative proportional fitting \citep{demingstephan}, read by \citet{csiszar75} as
alternating $I$-projections. Our windows also carry inequalities, whose projection
solves the Karush--Kuhn--Tucker system of Section~\ref{sec:assumptions}.

The preceding results hold a baseline $\mu$ fixed. They explain what conditioning
does within that component, but they do not ask which components survive when the
baseline itself is random. That is the question created by a de Finetti mixture.
Here the event changes the directing law through $\Pi_{N,E}$ as well as changing
each component law. Theorem~\ref{thm:master} controls the componentwise change for
a target anywhere in the horizon, and Theorem~\ref{thm:limit} identifies the random
directing measure that remains after mixing. Its projected law is $\PiE$, and
prediction under the limiting process is ordinary Bayesian prediction. We use
windows because an exact empirical equality can miss the type lattice altogether;
windows also accommodate inequalities without a separate construction.
\citet{jaynes}, \citet{zabell}, and, more recently,
\citet{foley2025} supply the maximum-entropy and Bayesian lineage. Closest to the
present setting, \citet{diaconisfreedman1988} prove conditional limit theorems for
exponential families and finite versions of de Finetti's theorem, with a sharp
$O(m/N)$ local expansion in the regular exact-sum setting (here $m$ is the block
size and $k$ is the alphabet size); their earlier paper \citep{diaconisfreedman}
gives the finite-$N$ de Finetti approximation. Predictive characterizations in
which parameters emerge as limits of predictive functionals are developed by
\citet{fortiniladelliregazzini} and reviewed by \citet{fortinipetrone}; their
primitive is predictive sufficiency and predictive-limit structure, whereas ours is
an active empirical constraint whose local component law is selected by Gibbs
conditioning.

The novelty is thus in the constrained-horizon formulation and its consequences
for the mixing law, not in exponential tilting itself.
Sections~\ref{sec:experiment}--\ref{sec:local} define the experiment and allow the
prediction target to occupy any positions in the active horizon.
Proposition~\ref{prop:split} and Theorem~\ref{thm:master} then split the constrained
type ensemble into two parts: a neighborhood of the $I$-projection handled by a
finite-exchangeability coupling, and a complement controlled by the
Csisz\'ar--Pinsker divergence gap. This gives an elementary finite-$N$ bound for
general polyhedral windows. Theorem~\ref{thm:limit} passes from those componentwise
bounds to the projected de Finetti law, whose directing measure is
$\Pstar_\mu$ for $\mu$ drawn from the limiting constraint law.
Corollary~\ref{cor:bayespred} supplies the resulting Bayesian predictive, and
Proposition~\ref{prop:dichotomy} separates the two possible limiting regimes.

Conditional exponential tilting, equivalence of ensembles, and finite de Finetti
consequences are classical, and the elementary global type argument is not sharper
than the classical regular conditional-limit results. A companion paper
\citep{polsonzantedeschi} takes the complementary view: it fixes a single baseline
and a finite measurable partition, then studies the local geometry of the collapse,
including Gaussian localization on the constraint tangent space, curvature
refinements, and estimation for parameter-indexed moment families. It carries no
outer prior over the baseline component $\mu$ and identifies no directing measure;
the present paper carries that prior and fixes no partition. The notation differs:
in that paper the horizon is
$n$, the window tolerance is $\eta_n$, the baseline is written $Q$, and a finite
measurable partition plays the role that the alphabet $\X$ plays here.

\begin{figure*}[t]
\centering
\includegraphics[width=\textwidth]{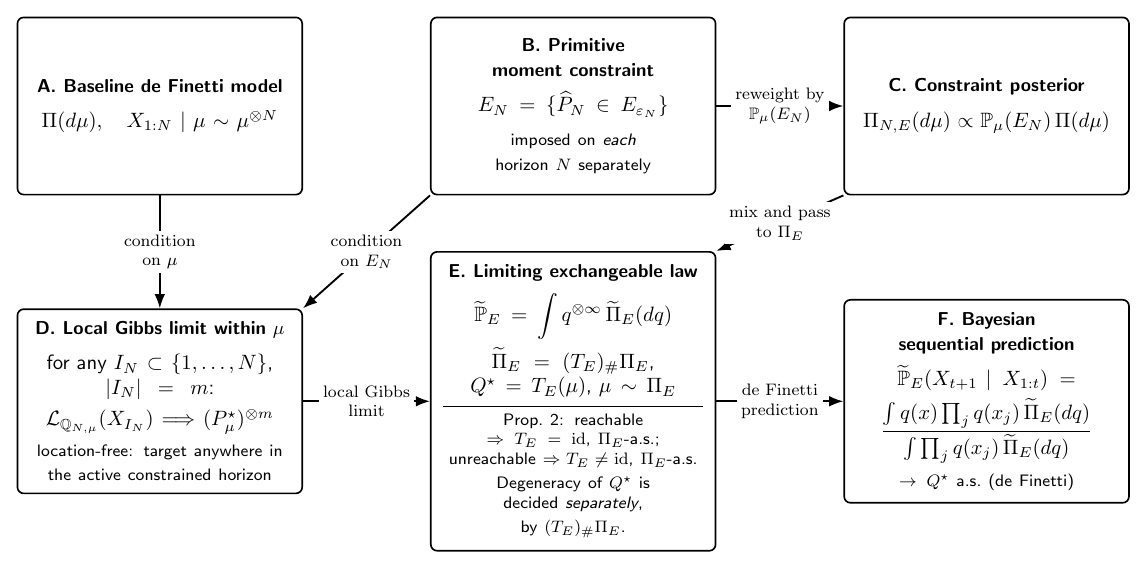}
\caption{The constraint is defined separately for each active finite horizon
$E_N=\{\widehat{P}_N\in E_{\varepsilon_N}\}$, giving the constrained-horizon law $\Qbb_N$ and
the constraint posterior $\Pi_{N,E}$. A prediction target may occupy arbitrary
coordinates in that horizon. Conditionally on $\mu$, Gibbs conditioning sends the
local law to the product $I$-projection $(\Pstar_\mu)^{\otimes m}$; the projected
component laws $\Pstar_\mu$ form the directing measures $\PiE$ of the limiting
exchangeable process, under which sequential prediction is Bayesian.}
\label{fig:fourcorner}
\end{figure*}

Section~\ref{sec:background} develops the de Finetti, type, and Gibbs-conditioning
ingredients in a common notation; Section~\ref{sec:master} proves the finite-horizon
master inequality, and Section~\ref{sec:limit} uses posterior-averaged component
control to obtain the projected exchangeable limit, with the master inequality
supplying a sufficient condition.
Section~\ref{sec:thermo} gives thermodynamic and local-fluctuation interpretations,
Section~\ref{sec:examples} develops four examples, and
Section~\ref{sec:discussion} summarizes the scope and open problems.
Figure~\ref{fig:fourcorner} summarizes the construction.

\section{Exchangeable prediction under moment constraints}
\label{sec:background}

Throughout, $\X=\{1,\dots,k\}$, $\simplex(\X)$ is the probability simplex, and for
$Q,P\in\simplex(\X)$ with $P>0$ the relative entropy is
$\KL{Q}{P}=\sum_x Q(x)\log\{Q(x)/P(x)\}$ and the Shannon entropy is
$H(Q)=-\sum_x Q(x)\log Q(x)$. We use total variation
$\TV{P-Q}=\tfrac12\sum_x|P(x)-Q(x)|$, and state every constant in this convention.

\subsection{De Finetti and the constraint posterior}
\label{sec:definetti}

\begin{theorem}[de Finetti; directing measure]
\label{thm:definetti}
Let $(X_i)_{i\ge1}$ be infinitely exchangeable in $\X$. There is a random measure
$\mu$ with law $\Pi$ such that, conditionally on $\mu$, the $X_i$ are i.i.d.\
$\mu$; moreover $\widehat{P}_N\to\mu$ a.s., and the conditional law of $\mu$ given the
observed information is the posterior.
\end{theorem}

The constrained-horizon law \eqref{eq:QN} carries the constraint posterior
$\Pi_{N,E}$ of \eqref{eq:constraintpost}, which we always denote $\Pi_{N,E}$;
$\Pi_E$ is reserved for a weak limit. We call $\Pi_{N,E}$ the constraint posterior
because the constraint is the primitive specification, not an observation made
once and discarded.

\subsection{The method of types}
\label{sec:types}

A \emph{type} is an empirical measure achievable with denominator $N$,
$\TN=\{n/N:n\in\mathbb{N}_0^k,\ \sum_j n_j=N\}$. There are at most
$(N+1)^k$ types \citep[Thm.~11.1.1]{coverthomas}. For i.i.d.\ $P>0$
\citep[Thm.~11.1.4]{coverthomas},
\begin{equation}
\label{eq:typeprob}
(N+1)^{-k}\,e^{-N\KL{Q}{P}}\le \Pr_P(\widehat{P}_N=Q)\le e^{-N\KL{Q}{P}}
\end{equation}
for every type $Q\in\TN$. The bound \eqref{eq:typeprob} is a counting bound: there
are many types, and their probabilities must be controlled
uniformly. Section~\ref{sec:hs} records the case in which no counting is needed, and
where the method of types is visible in its most elementary form.

\subsection{\texorpdfstring{$I$}{I}-projections and Gibbs conditioning}
\label{sec:iproj}

For closed convex $E\subseteq\simplex(\X)$ and $P>0$, the $I$-projection is
\begin{equation}
\label{eq:iproj}
\Pstar=\argmin_{Q\in E}\KL{Q}{P},
\end{equation}
unique by strict convexity of $\KL{\cdot}{P}$; for an affine moment set with moment
functions $g$ (Section~\ref{sec:windows}) and an interior projection, it is the
finite exponential tilt
$\Pstar(x)\propto P(x)e^{\lambda^\top g(x)}$ \citep{csiszar75}; boundary projections
are obtained as limiting tilts.
The Pythagorean inequality of \citet{csiszar75}, for convex $E$ and $Q\in E$,
$\KL{Q}{P}\ge\KL{Q}{\Pstar}+\KL{\Pstar}{P}$, combined with Pinsker's inequality
$\KL{Q}{\Pstar}\ge 2\TV{Q-\Pstar}^2$, gives the divergence-gap bound
\begin{equation}
\label{eq:klgap}
\KL{Q}{P}-\KL{\Pstar}{P}\ge 2\TV{Q-\Pstar}^2,\qquad Q\in E .
\end{equation}
The conditional (Gibbs) form of Sanov's theorem \citep{csiszar84,leonardnajim}
states that, under the usual conditioning-set regularity and along horizons for
which $\{\widehat{P}_N\in E\}$ has positive probability (or for an appropriate
neighborhood formulation), the conditional law of a fixed block converges in total
variation to $(\Pstar)^{\otimes m}$ when $E$ is convex and $\Pstar$ is unique.
Section~\ref{sec:master} makes this effective and location-free for shrinking
polyhedral windows. The projection $\Pstar$ is a
maximum-\emph{entropy} law only when $P$ is uniform; in general it is the
minimum-relative-entropy, or $I$-projection, law.

\subsection{Type decomposition and the Heath--Sudderth construction}
\label{sec:decomp}

For a fixed block size $m\le N$, writing $\Pr(x_{1:m}\mid\cdot)$ for
$\Pr(X_{1:m}=x_{1:m}\mid\cdot)$, the law of total probability over types gives
\begin{equation}
\label{eq:hs}
\Pr(x_{1:m})=\sum_{Q\in\TN}\Pr(x_{1:m}\mid \widehat{P}_N=Q)\,\Pr(\widehat{P}_N=Q).
\end{equation}
Conditionally on $\widehat{P}_N=Q$, the block is a draw without replacement from a
population of composition $Q$, and
\begin{equation}
\label{eq:coupling}
\TV{\Pr(X_{1:m}\in\cdot\mid \widehat{P}_N=Q)-Q^{\otimes m}}\le \frac{m(m-1)}{2N}.
\end{equation}
This is an elementary collision coupling. Draw $m$ indices from the population of $N$
uniformly \emph{with} replacement; conditionally on the event that no index repeats, the
sample is distributed as a draw \emph{without} replacement. The two index vectors
therefore differ in total variation by at most the collision probability, at most
$\binom{m}{2}/N=m(m-1)/(2N)$, and the labels are a deterministic function of the
indices, so the same bound holds for them. The with-replacement labels are i.i.d.\ $Q$,
which is \eqref{eq:coupling}. \citet{diaconisfreedman} prove the general
finitely-exchangeable statement, that a finitely exchangeable law is close in total
variation to a mixture of product measures; \eqref{eq:coupling} is the within-type
special case.

Equation \eqref{eq:hs} is exact and structural, and \eqref{eq:coupling} says that a
fixed block, conditioned on the empirical measure, is an independent sample from that
empirical measure up to $O(m^2/N)$. Together, these two identities underlie the
analysis: within a fixed horizon the type decomposition supplies the mixture
structure, Gibbs conditioning selects the relevant components, and the master
inequality of Section~\ref{sec:master} quantifies their combination.

\subsubsection{The method of types on a scalar lattice}
\label{sec:hs}

Before imposing the constraint, it is useful to see what \eqref{eq:hs} and
\eqref{eq:coupling} already deliver on their own. \citet{heathsudderth} use them to
prove de Finetti's theorem for binary sequences, and their argument is the method of
types in its most transparent form. We reproduce it because it is the skeleton of
Theorem~\ref{thm:limit}.

Take $k=2$. A type is then the scalar $\widehat{P}_N=S_N/N$ with
$S_N=\sum_{i\le N}X_i$, the lattice $\TN$ is $\{0,1/N,\dots,1\}$, and the law of the
type is summarized by the distribution function
$F_N(y)=\Pr(\widehat{P}_N\le y)$, $y\in[0,1]$, a right-continuous step function with at
most $N+1$ jumps. Write $s=\sum_{i\le m}x_i$. Exchangeability makes the conditional
block law hypergeometric, a draw without replacement from an urn with $j$ ones and
$N-j$ zeros, and \eqref{eq:hs} becomes a Stieltjes integral against $F_N$:
\begin{equation}
\label{eq:hsstieltjes}
\Pr(X_{1:m}=x_{1:m})=\int_0^1\varphi_{N,s}(y)\,dF_N(y),
\end{equation}
with $\varphi_{N,s}(y)=(Ny)_s\,(N-Ny)_{m-s}/(N)_m$ and
$(a)_b=a(a-1)\cdots(a-b+1)$. The identity holds for every $N\ge m$, and its left side
does not depend on $N$. Evaluated at a single configuration, \eqref{eq:coupling} says
that the integrand converges to a polynomial uniformly on the lattice,
\begin{equation}
\label{eq:hsuniform}
\sup_{y\in\TN}\big|\varphi_{N,s}(y)-y^s(1-y)^{m-s}\big|\le\frac{m(m-1)}{2N}.
\end{equation}

\begin{theorem}[de Finetti for binary sequences; Heath and Sudderth]
\label{thm:hs}
Let $(X_i)_{i\ge1}$ be infinitely exchangeable with values in $\{0,1\}$. There is a
unique distribution function $F$ on $[0,1]$ with
\begin{equation}
\label{eq:hsdefinetti}
\Pr(X_1=x_1,\dots,X_m=x_m)=\int_0^1 y^s(1-y)^{m-s}\,dF(y),
\end{equation}
$s=\sum_{i\le m}x_i$, for every $m$ and every $x_{1:m}$, and $F_N\Rightarrow F$, the
law of the directing measure.
\end{theorem}

\begin{proof}
The type lattice sits in the compact interval $[0,1]$, so by the Helly selection
theorem \citep[Theorem~25.9]{billingsley} there is a subsequence along which
$F_{N_j}\Rightarrow F$ for a proper distribution function $F$, no mass escaping the
compact interval. Fix $x_{1:m}$ and set $g_s(y)=y^s(1-y)^{m-s}$, continuous and bounded
on $[0,1]$. By \eqref{eq:hsstieltjes} and the triangle inequality,
$|\Pr(X_{1:m}=x_{1:m})-\int_0^1 g_s\,dF|$ is at most
$\sup_{y\in\mathcal{T}_{N_j}}|\varphi_{N_j,s}(y)-g_s(y)|$, which is
$\le m(m-1)/(2N_j)$ by \eqref{eq:hsuniform} and suffices because $F_{N_j}$ is carried
by the lattice, plus $|\int_0^1 g_s\,dF_{N_j}-\int_0^1 g_s\,dF|$, which tends to zero by
the Helly--Bray theorem \citep[Theorem~25.8(ii)]{billingsley}. Both terms vanish along
the subsequence while the left side does not depend on $j$, so it is zero, which is
\eqref{eq:hsdefinetti}. Taking $x_{1:m}=(1,\dots,1)$ gives $\int_0^1 y^m dF(y)=
\Pr(X_1=\dots=X_m=1)$ for every $m$, so the moments of $F$ are determined by the law
of the sequence; a probability measure on a compact interval is determined by its
moments, which gives uniqueness of $F$ and hence convergence of the full sequence.
\end{proof}

\subsubsection{Counting and weighting of types}
\label{sec:counting}

The binary argument isolates the part of the method of types that does not require
counting. Because a binary type is a scalar in $[0,1]$, its weights are encoded by the
single monotone function $F_N$, and the exponential estimate \eqref{eq:typeprob} is
unnecessary. Compactness supplies a weakly convergent subsequence;
\eqref{eq:coupling} replaces the hypergeometric integrand uniformly by its binomial
counterpart; and Helly--Bray passes to the limit. The directing distribution thereby
arises as the weak limit of the empirical-mean laws, rather than being postulated and
subsequently identified.

Conditioning is where the weights require separate attention. In the
unconstrained identity \eqref{eq:hs}, the probabilities
$\Pr(\widehat{P}_N=Q)$ have already been absorbed into $F_N$. After imposing
$E_N$, the relevant weights are
\[
w_N(Q)=\Pr(\widehat{P}_N=Q\mid E_N),
\]
and their denominator may itself be small. Lemma~\ref{lem:gap} controls that
denominator and the numerator together. It combines the type bound
\eqref{eq:typeprob} with the Csisz\'ar--Pinsker gap \eqref{eq:klgap}, showing that
types separated from the $I$-projection carry exponentially little conditional
mass. In this sense the split-set argument of Section~\ref{sec:master} is a
constrained Heath--Sudderth decomposition: it exposes the weights and recenters
the component law at the projection. Theorem~\ref{thm:limit} uses precisely this
finite-$N$ control.

\section{The split-set proof and master inequality}
\label{sec:master}

Work conditionally on the directing measure $\mu=P$, so $X_{1:N}$ is i.i.d.\ $P$.

\subsection{Moment windows}
\label{sec:windows}

Let $g_1,\dots,g_r:\X\to\Real$ and thresholds $\alpha_j$ define the polyhedral set
\begin{equation}
\label{eq:momentset}
\begin{aligned}
E=\big\{Q\in\simplex(\X):\ &\langle g_j,Q\rangle=\alpha_j\ (j\le r_0),\\[-1pt]
&\langle g_j,Q\rangle\ge\alpha_j\ (r_0<j\le r)\big\},
\end{aligned}
\end{equation}
with $\langle g,Q\rangle=\sum_x g(x)Q(x)$, and the window
\begin{equation}
\label{eq:window}
\begin{aligned}
E_\varepsilon=\big\{Q:\ &|\langle g_j,Q\rangle-\alpha_j|\le\varepsilon\ (j\le r_0),\\[-1pt]
&\langle g_j,Q\rangle\ge\alpha_j-\varepsilon\ (r_0<j\le r)\big\}.
\end{aligned}
\end{equation}
The only lattice issue is feasibility. Rounding a probability vector to an $N$-type
moves each coordinate by $O(N^{-1})$ and hence changes every moment
$\langle g_j,\cdot\rangle$ by at most $c_{\mathrm{rnd}}/N$, where
$c_{\mathrm{rnd}}=k\max_j\|g_j\|_\infty$. Thus a point in the inner window
$E_{\varepsilon-c_{\mathrm{rnd}}/N}$ can be rounded to an $N$-type in
$E_\varepsilon$. We use this observation in the denominator estimate below.

\subsection{Standing assumptions}
\label{sec:assumptions}

\begin{description}
\item[(A1)] \emph{Finite alphabet.} $\X=\{1,\dots,k\}$, with $k<\infty$.
\item[(A2)] \emph{Positive baseline.} $\min_x P(x)>0$.
\item[(A3)] \emph{Polyhedral constraint.} The set $E$ in
\eqref{eq:momentset} is nonempty, closed, and convex.
\item[(A4)] \emph{Interior projection.} The $I$-projection $\Pstar$ in
\eqref{eq:iproj} satisfies $\min_x\Pstar(x)>0$. Its uniqueness follows from the
strict convexity of $\KL{\cdot}{P}$.
\item[(A5)] \emph{Window stability.} If $\Pstar_\varepsilon$ is the
$I$-projection of $P$ onto $E_\varepsilon$, then, for some $C_0<\infty$,
$\TV{\Pstar_\varepsilon-\Pstar_{\varepsilon'}}\le C_0|\varepsilon-\varepsilon'|$
whenever $\varepsilon,\varepsilon'\ge0$ are sufficiently small.
\end{description}

We use (A5) only through this Lipschitz bound, which follows from standard
strong-regularity conditions for the locally active parametric
Karush--Kuhn--Tucker system \citep{robinson}. In particular,
$\TV{\Pstar_{\varepsilon_N-c_{\mathrm{rnd}}/N}-\Pstar_{\varepsilon_N}}=O(1/N)$.
Interiority (A4) and positivity (A2) make $\KL{\cdot}{P}$ locally Lipschitz near
$\Pstar$, which the denominator argument uses.

\subsection{The split-set bound and the master inequality}
\label{sec:masterthm}

The proof uses a two-region decomposition of the constrained type space. Types
within total-variation radius $\delta$ of the window $I$-projection
$P^\bullet=\Pstar_{\varepsilon_N}$ are controlled geometrically: their conditional
block laws are close to the corresponding product law, up to the finite-population
sampling error. Types outside that radius are controlled probabilistically: the
Csisz\'ar--Pinsker gap \eqref{eq:klgap} imposes an exponential relative-entropy
penalty. We call this the \emph{split-set argument}, and write the admissible
types as
\begin{equation}
\label{eq:splitset}
\begin{aligned}
G_{N,\delta}&=\{Q\in E_{\varepsilon_N}\cap\TN:\TV{Q-P^\bullet}\le\delta\},\\
B_{N,\delta}&=\{Q\in E_{\varepsilon_N}\cap\TN:\TV{Q-P^\bullet}>\delta\}.
\end{aligned}
\end{equation}
For an $m$-element index set $I_N\subseteq\{1,\dots,N\}$,
$L_{N,I_N,P}=\mathcal{L}_P(X_{I_N}\mid E_N)$ is the constrained local law; with
$L_Q=\mathcal{L}_P(X_{1:m}\mid \widehat{P}_N=Q)$ and $w_N(Q)=\Pr_P(\widehat{P}_N=Q\mid E_N)$ it splits
into the subprobability contributions
\begin{equation}
\label{eq:GBcontrib}
\begin{aligned}
L^{G}_{N,I_N,P}(A)&=\sum_{Q\in G_{N,\delta}}w_N(Q)L_Q(A),\\
L^{B}_{N,I_N,P}(A)&=\sum_{Q\in B_{N,\delta}}w_N(Q)L_Q(A),
\end{aligned}
\end{equation}
of masses $1-\beta_N$ and $\beta_N=\Pr_P(\widehat{P}_N\in B_{N,\delta}\mid E_N)$, with
$L_{N,I_N,P}=L^{G}_{N,I_N,P}+L^{B}_{N,I_N,P}$.

\begin{proposition}[Split-set decomposition]
\label{prop:split}
Under \textup{(A1)--(A5)} with $N\varepsilon_N>c_{\mathrm{rnd}}$ eventually, for
all sufficiently large $N$, every $\delta>0$, and every $m$-element $I_N$,
\begin{equation}
\label{eq:splitbound}
\TV{L_{N,I_N,P}-(P^\bullet)^{\otimes m}}\le m\delta+\frac{m(m-1)}{2N}+\beta_N,
\end{equation}
and the bad-set mass obeys
$\beta_N=\Pr_P(\widehat{P}_N\in B_{N,\delta}\mid E_N)\le C(P,E)(N+1)^{2k}e^{-2N\delta^2}$.
\end{proposition}

Adding the window-stability term, which bounds
$\TV{(P^\bullet)^{\otimes m}-(\Pstar)^{\otimes m}}$ by $C_0 m\varepsilon_N$, turns
Proposition~\ref{prop:split} into the master inequality.

\needspace{15\baselineskip}
\begin{theorem}[Componentwise constrained-horizon master inequality]
\label{thm:master}
Assume \textup{(A1)--(A5)} and $\varepsilon_N\downarrow0$ with
$N\varepsilon_N>c_{\mathrm{rnd}}$ eventually. There is $C(P,E)<\infty$ such that,
for all sufficiently large $N$, every $m$-element $I_N\subseteq\{1,\dots,N\}$, and
every $\delta>0$,
\begin{equation}
\label{eq:master}
\begin{aligned}
\TV{L_{N,I_N,P}&-(\Pstar)^{\otimes m}}\\
\le{}& \underbrace{m\delta+\tfrac{m(m-1)}{2N}}_{\text{good-set approximation}}\\
&{}+\underbrace{C(P,E)(N+1)^{2k}e^{-2N\delta^2}}_{\text{bad-set mass}}\\
&{}+\underbrace{C_0\,m\varepsilon_N}_{\text{window transfer}}.
\end{aligned}
\end{equation}
\end{theorem}

In \eqref{eq:master}, $G_{N,\delta}$ contributes the local approximation and the
sampling-without-replacement error, $B_{N,\delta}$ the exponentially small
conditional mass, and the passage $P^\bullet\to\Pstar$ the window-transfer error.
The bound is independent of the positions in $I_N$: by \eqref{eq:locfree} the
constrained component law is exchangeable, so $I_N=\{1,\dots,m\}$ gives the general
case exactly. The proof is in Appendix~\ref{app:proofs}.

\begin{corollary}[A tuning]
\label{cor:tuning}
Under the hypotheses of Theorem~\ref{thm:master}, taking
$\delta_N=\sqrt{(k+1)\log(N+1)/N}$ gives, for fixed $m$,
\begin{equation}
\label{eq:rate}
\TV{L_{N,I_N,P}-(\Pstar)^{\otimes m}}
=O\!\Big(m\varepsilon_N+m\sqrt{\tfrac{\log N}{N}}+\tfrac{m^2}{N}\Big).
\end{equation}
\end{corollary}

The choice $\delta_N$ is one tuning that makes the exponential remainder
$O((N+1)^{-2})$. In the regular exact-sum
setting \citet{diaconisfreedman1988} obtain a sharper $O(m/N)$ expansion; the
present elementary argument trades sharpness for the generality of polyhedral
windows. The $\sqrt{\log N/N}$ term is an artifact of the union bound over the at
most $(N+1)^k$ types and is absent in the regular exact-sum settings of
\citet{diaconisfreedman1988}, where the leading term
carries no logarithm.

Theorem~\ref{thm:master} is a statement conditional on the baseline $\mu$, and it is
the layer that Theorem~\ref{thm:limit} then mixes. The same estimate, in the
partition-chart setting with chart-dependent constants and no prior, is the
window-robust collapse bound of the companion paper \citep{polsonzantedeschi}, where
it serves a different end: there it licenses an induced conditional likelihood on the
chart, here it licenses the passage to a limiting exchangeable law.

\section{The projected de Finetti limit}
\label{sec:limit}

We now integrate over the directing measure. The exact finite-$N$ identity is
\begin{equation}
\label{eq:exactmix}
\Qbb_N(X_{I_N}\in\cdot)=\int \Qbb_{N,\mu}(X_{I_N}\in\cdot)\,\Pi_{N,E}(d\mu).
\end{equation}
Write $T_E(\mu)=\Pstar_\mu=\argmin_{Q\in E}\KL{Q}{\mu}$ for the projection map.

\begin{theorem}[Projected de Finetti limit]
\label{thm:limit}
Suppose \textup{(i)} $\Pi_{N,E}\Rightarrow\Pi_E$ weakly; \textup{(ii)} $T_E$ is
continuous on a neighborhood of $\operatorname{supp}(\Pi_E)$; and \textup{(iii)}
for every fixed $m$,
\begin{equation}
\label{eq:posavg}
\int\TV{\Qbb_{N,\mu}(X_{1:m}\in\cdot)-T_E(\mu)^{\otimes m}}\,\Pi_{N,E}(d\mu)\to0.
\end{equation}
Then for every fixed $m$ and every sequence of $m$-element sets
$I_N\subseteq\{1,\dots,N\}$,
\begin{equation}
\label{eq:mixlimit}
\Qbb_N(X_{I_N}\in\cdot)\ \longrightarrow\ \int T_E(\mu)^{\otimes m}\,\Pi_E(d\mu)
\end{equation}
in total variation. Consequently, with $\PiE=(T_E)_\#\Pi_E$, the family
$\widetilde{\Pr}^{(m)}_E=\int q^{\otimes m}\PiE(dq)$ is projectively consistent and
defines an infinite exchangeable law
$\widetilde{\Pr}_E=\int q^{\otimes\infty}\PiE(dq)$, whose random directing measure
is $Q^\star=T_E(\mu)=\Pstar_\mu$, $\mu\sim\Pi_E$.
\end{theorem}

\begin{proof}
Fix $x_{1:m}$. By \eqref{eq:exactmix} and \eqref{eq:locfree}, the finite-horizon
marginal is the corresponding mixture over $\Pi_{N,E}$. The triangle inequality
gives
\[
\begin{aligned}
&\TV{\Qbb_N(X_{I_N}\in\cdot)-\textstyle\int T_E(\mu)^{\otimes m}\Pi_E(d\mu)}\\
&\ \le \int\TV{\Qbb_{N,\mu}(X_{1:m}\in\cdot)-T_E(\mu)^{\otimes m}}\Pi_{N,E}(d\mu)\\
&\quad{}+\TV{\textstyle\int T_E(\mu)^{\otimes m}\Pi_{N,E}(d\mu)-\int T_E(\mu)^{\otimes m}\Pi_E(d\mu)}.
\end{aligned}
\]
The first term vanishes by \eqref{eq:posavg}. In the second, the function
$f_{x_{1:m}}(\mu)=\prod_{i=1}^m\Pstar_\mu(x_i)$ is bounded by $1$ and continuous at
every point of $\operatorname{supp}(\Pi_E)$ by (ii), so its discontinuity set has
$\Pi_E$-measure zero and the bounded almost-everywhere-continuous form of the
Portmanteau theorem, with weak convergence (i), gives
$\int f_{x_{1:m}}\,d\Pi_{N,E}\to\int f_{x_{1:m}}\,d\Pi_E$; summing over the finitely
many $x_{1:m}$ gives convergence in total variation. The limiting family
$\widetilde{\Pr}^{(m)}_E$ is a de Finetti mixture of i.i.d.\ laws, hence
projectively consistent, and Kolmogorov extension gives $\widetilde{\Pr}_E$ with
directing measure $Q^\star$.

This proof has the same structure as Theorem~\ref{thm:hs}, with the projection
inserted. There \eqref{eq:hsstieltjes} is the finite-$N$ mixture representation,
$F_N$ is the mixing law, finite exchangeability controls the within-component
error, and Helly--Bray handles the passage across components. Here these roles are
played by \eqref{eq:exactmix}, the constraint posterior $\Pi_{N,E}$, hypothesis
(iii) as supplied by Theorem~\ref{thm:master}, and the Portmanteau theorem applied
to $f_{x_{1:m}}$. Gibbs conditioning changes the component law from $\mu$ to
$\Pstar_\mu$ but preserves the mixture argument. The additional price is continuity
of $T_E$ in (ii), a condition that is vacuous in the untilted case because the
integrand is then a polynomial.
\end{proof}

The finite-dimensional marginal $\int(\Pstar_\mu)^{\otimes m}\Pi_E(d\mu)$ is
deterministic. Under the limiting law, the random object is the directing measure
$Q^\star=\Pstar_\mu$, with $\mu\sim\Pi_E$ and hence $Q^\star\sim\PiE$, as required
by exchangeability.

\begin{corollary}[Uniform sufficient condition]
\label{cor:uniform}
In addition to \textup{(i)--(ii)} of Theorem~\ref{thm:limit}, suppose there is a
compact $K\subseteq\simplex(\X)$ with $\Pi_{N,E}(K^c)\to0$,
$\inf_{\mu\in K,x}\mu(x)>0$, the projections $\Pstar_\mu$ uniformly interior on
$K$, and the constants of Theorem~\ref{thm:master} uniform over $K$. Then
$\sup_{\mu\in K}\TV{\Qbb_{N,\mu}(X_{1:m}\in\cdot)-T_E(\mu)^{\otimes m}}\to0$, the
$K^c$ contribution to \eqref{eq:posavg} is at most $\Pi_{N,E}(K^c)\to0$, so
\eqref{eq:posavg} holds. The weak-convergence and continuity hypotheses
\textup{(i)--(ii)} are still required.
\end{corollary}

\subsection{Bayesian prediction under the constrained limit}
\label{sec:bayespred}

\begin{corollary}[Bayesian prediction under the constrained limit]
\label{cor:bayespred}
Under $\widetilde{\Pr}_E$, for fixed $t$ and a history $x_{1:t}$ of positive
marginal probability,
\begin{equation}
\label{eq:bayespred}
\begin{aligned}
&\widetilde{\Pr}_E(X_{t+1}=x\mid X_{1:t}=x_{1:t})\\
&\qquad=\frac{\int \Pstar_\mu(x)\prod_{j=1}^t\Pstar_\mu(x_j)\,\Pi_E(d\mu)}
{\int\prod_{j=1}^t\Pstar_\mu(x_j)\,\Pi_E(d\mu)} .
\end{aligned}
\end{equation}
This is the standard de Finetti posterior predictive under the induced prior
$\PiE=(T_E)_\#\Pi_E$, and
$\widetilde{\Pr}_E(X_{t+1}\in\cdot\mid X_{1:t})\to Q^\star$ almost surely, where
$Q^\star$ has law $\PiE$.
\end{corollary}

Corollary~\ref{cor:bayespred} is the cleanest formal expression of the title:
Sanov determines the transformed directing component $\Pstar_\mu$, and de Finetti
supplies Bayesian prediction under the resulting exchangeable law.

\subsection{Reachability and degeneracy of the limit}
\label{sec:dichotomy}

Theorem~\ref{thm:limit} identifies the limit in terms of $\Pi_E$ and $T_E$. Both are
determined by the prior, and the outcome is a dichotomy that governs every example
below. The constraint posterior is exponentially weighted: for each fixed strictly
positive $\mu$, under the window and type conditions of Section~\ref{sec:master},
\begin{equation}
\label{eq:constraintLDP}
\begin{aligned}
\tfrac1N\log\Pr_\mu(E_N)&\ \longrightarrow\ -I_E(\mu),\\
\Pi_{N,E}(d\mu)&\propto e^{-N I_E(\mu)+o_\mu(N)}\,\Pi(d\mu),
\end{aligned}
\end{equation}
pointwise in $\mu$, with $I_E(\mu)=\inf_{Q\in E}\KL{Q}{\mu}$.

\begin{proposition}[Reachability and degeneracy]
\label{prop:dichotomy}
Assume hypothesis \textup{(i)} of Theorem~\ref{thm:limit}. Assume also that
$\operatorname{supp}\Pi\subseteq\{\mu:\min_x\mu(x)>0\}$, which is the regime in which
\eqref{eq:constraintLDP} was derived and which gives $I^\star<\infty$, and that the
Laplace principle behind \eqref{eq:constraintLDP} holds \emph{uniformly} on
$\operatorname{supp}\Pi$, that is, $o_\mu(N)=o(N)$ uniformly there. Write
$I^\star=\inf_{\operatorname{supp}\Pi}I_E$, and note that $I_E$ vanishes exactly on
$E$, since $\KL{Q}{\mu}=0$ only at $Q=\mu$.
\begin{enumerate}[label=\textup{(\roman*)}]
\item \emph{Support.} $\Pi_E$ is carried by $\operatorname{argmin}_{\operatorname{supp}\Pi}I_E
=\{I_E=I^\star\}$.
\item \emph{Reachable.} If $I^\star=0$, then $\Pi_E$ is carried by $E$ and
$T_E=\mathrm{id}$ $\Pi_E$-almost surely, so $\PiE=\Pi_E$: the constraint restricts the
prior and the $I$-projection is asymptotically inactive.
\item \emph{Unreachable.} If $I^\star>0$, the surviving components are genuine
projections, $\Pstar_\mu\neq\mu$ for $\Pi_E$-almost every $\mu$.
\item \emph{Degeneracy.} Let $E$ be affine, let the components of
$\operatorname{supp}\Pi_E$ be strictly positive with interior $I$-projections, in the
sense of \textup{(A2)} and \textup{(A4)} applied componentwise, and let $g_1,\dots,g_{r_0}$ have rank $r_0$
\emph{modulo constants}. Then $\Pstar_\mu=\Pstar_\nu$ if and only if $\mu$ and $\nu$ lie in the same
orbit of the exponential family generated by $g_1,\dots,g_{r_0}$. Hence $\PiE$ is
degenerate if and only if $\operatorname{supp}\Pi_E$ lies in a single tilt orbit;
otherwise the surviving randomness is the component of $\Pi_E$ transverse to the tilt
family, the orbit space having dimension $k-1-r_0$.
\end{enumerate}
\end{proposition}

\begin{proof}
(i) The divergence $\KL{\cdot}{\cdot}$ is jointly lower semicontinuous on
$\simplex(\X)$ and $E$ is compact, so $I_E$ is lower semicontinuous and
$U_\delta=\{I_E>I^\star+\delta\}$ is open for every $\delta>0$.
As a closed subset of the compact simplex,
$\operatorname{supp}\Pi$ is itself compact, so lower semicontinuity also gives a
minimizer: there is $\mu_0\in\operatorname{supp}\Pi$ with $I_E(\mu_0)=I^\star$. Fix $\delta>0$. Let $Q_0=\Pstar_{\mu_0}\in E$ attain $\KL{Q_0}{\mu_0}=I^\star$, finite
because $\operatorname{supp}\Pi$ consists of strictly positive components, so
$\mu_0(x)>0$ wherever $Q_0(x)>0$. Then $\mu\mapsto\KL{Q_0}{\mu}$ is continuous at
$\mu_0$ and $I_E\le\KL{Q_0}{\cdot}$ pointwise, so $I_E<I^\star+\delta/2$ on a
neighborhood $G$ of $\mu_0$, with $\Pi(G)>0$ because
$\mu_0\in\operatorname{supp}\Pi$. By \eqref{eq:constraintLDP} with uniform $o(N)$,
\[
\begin{aligned}
\frac{\Pi_{N,E}(U_\delta)}{\Pi_{N,E}(G)}
&\le\frac{e^{-N(I^\star+\delta)+o(N)}}
{\Pi(G)\,e^{-N(I^\star+\delta/2)+o(N)}}\\
&=e^{-N\delta/2+o(N)}\to0,
\end{aligned}
\]
so $\Pi_{N,E}(U_\delta)\to0$. As $U_\delta$ is open, the portmanteau theorem and
hypothesis (i) of Theorem~\ref{thm:limit} give
$\Pi_E(U_\delta)\le\liminf_N\Pi_{N,E}(U_\delta)=0$. Since
$\{I_E>I^\star\}=\bigcup_{n\ge1}U_{1/n}$ is a countable union, $\Pi_E(I_E>I^\star)=0$.
Finally $\Pi_E$ is carried by the closed set $\operatorname{supp}\Pi$, being a weak limit
of measures carried by it, so $\Pi_E$ is carried by
$\{I_E=I^\star\}\cap\operatorname{supp}\Pi=\operatorname{argmin}_{\operatorname{supp}\Pi}I_E$. (ii) If $I^\star=0$ then
$\{I_E=0\}=E\cap\operatorname{supp}\Pi$, and for $\mu\in E$ the minimizer of
$\KL{\cdot}{\mu}$ over $E$ is $\mu$ itself, so $T_E|_E=\mathrm{id}$. (iii) is immediate
from (i), since $I_E(\mu)>0$ forces $\Pstar_\mu\neq\mu$. (iv) For affine $E$,
$\Pstar_\mu(x)\propto\mu(x)e^{\lambda^\top g(x)}$. If $\nu\propto\mu\,e^{\theta^\top g}$
then $\Pstar_\nu\propto\mu\,e^{(\theta+\lambda')^\top g}$, and the multiplier meeting the
constraints returns $\Pstar_\nu=\Pstar_\mu$; conversely $\Pstar_\mu=\Pstar_\nu$ forces
$\mu$ and $\nu$ to differ by a tilt. Hence $T_E$ is constant exactly on orbits, $E$ is a
transversal, and the orbit space has dimension $(k-1)-r_0$. On the boundary two components can share a
projection without lying in a common finite tilt orbit; the rank is taken modulo
constants because adding a constant to $g_j$ is absorbed by the normalization.
\end{proof}

Proposition~\ref{prop:dichotomy} locates $\Pi_E$ and settles whether the projection
acts. It does not identify how $\Pi_E$ distributes its mass within
$\{I_E=I^\star\}$; that weighting depends on the window. Two identifications are
available, each under its own hypotheses and each proved directly from
\eqref{eq:constraintpost}, so neither inherits the support-positivity assumption of
Proposition~\ref{prop:dichotomy}. The examples of Section~\ref{sec:examples} use them
in that form: $\mathrm{Beta}(2,2)$ and $\mathrm{Dirichlet}(1,1,1)$ have topological
support meeting the boundary of the simplex, which they charge with probability zero.

\begin{corollary}[Truncation]
\label{cor:trunc}
Under the constraint posterior \eqref{eq:constraintpost}, assume the acceptance
condition
\begin{equation}
\label{eq:accept}
\begin{aligned}
\Pr_\mu(E_N)&\to1 && \text{$\Pi$-a.e.\ on }E,\\
\Pr_\mu(E_N)&\to0 && \text{$\Pi$-a.e.\ off }E .
\end{aligned}
\end{equation}
If $\Pi(E)>0$ then $\Pi_E=\Pi(\cdot\mid E)$ and $\PiE=\Pi(\cdot\mid E)$: the constraint
truncates the prior. Condition \eqref{eq:accept} holds when $\Pi(\partial E)=0$, the
interior and boundary taken in $\simplex(\X)$ and not relative to
$\operatorname{supp}\Pi$: for $\mu$ interior to $E$ the law of large numbers puts
$\PN$ in $E$ eventually, and for $\mu$ off the closed set $E$ it puts $\PN$ outside
$E_{\varepsilon_N}$ eventually. The inequality constraint of
Section~\ref{sec:bernoulli} satisfies it. An equality constraint does not, $E$ then being
its own boundary; for regular components with nondegenerate moment covariance
\eqref{eq:accept} is instead ensured by a window with
$\sqrt N\,\varepsilon_N\to\infty$, and it fails at the $N^{-1/2}$ window scale. A
degenerate component, with $\langle g,\cdot\rangle$ equal to $\alpha$ $\mu$-almost
surely, has acceptance one and satisfies \eqref{eq:accept} at any window.
\end{corollary}

\begin{proof}
Bounded convergence in \eqref{eq:constraintpost}. For the sufficiency claim, if
$\mu$ is interior to $E$ then the empirical moments $\langle g_j,\PN\rangle$ satisfy the
constraint strictly and a law of large numbers gives $\Pr_\mu(E_N)\to1$; off $E$ the
same argument gives $0$.
\end{proof}

\begin{corollary}[Co-area limit]
\label{cor:coarea}
Suppose $\Pi(E)=0$ and that $E=\{Q:\langle g,Q\rangle=\alpha\}$ is a full-rank
equality constraint. Disintegrate the prior along the moment map
$\mu\mapsto\langle g,\mu\rangle$,
\begin{equation}
\label{eq:disint}
\Pi(d\mu)=\Pi_a(d\mu)\,h(a)\,da\quad\text{locally near }a=\alpha,
\end{equation}
where $h$ is the density of the pushforward of $\Pi$ under that map and $\Pi_a$ is the
conditional law on the fiber $\{\mu:\langle g,\mu\rangle=a\}$. Assume $h$ is continuous at $\alpha$ with $h(\alpha)>0$, which
also makes the constraint reachable. Assume that a version of the disintegration
$a\mapsto\Pi_a$ is weakly continuous at $\alpha$, and fix that version; explicitly,
for every bounded continuous $f$ on $\simplex(\X)$,
$a\mapsto\int f(\mu)\,\Pi_a(d\mu)$ is continuous at $\alpha$. Thus $\Pi_\alpha$
denotes this continuous version on the prior-null fiber. In smooth-density settings
the canonical co-area disintegration is one such version. Assume also that
$\Omega(\mu)=\operatorname{Cov}_\mu(g)$ is continuous and uniformly nondegenerate near
$E$, and that the local central limit theorem for $\langle g,\widehat{P}_N\rangle$ holds
uniformly there \emph{at the normalizing scale}. Precisely, with
$\Delta(\mu)=\langle g,\mu\rangle-\alpha$, $s(\mu)=\{\Omega(\mu)/N\}^{1/2}$ and
\begin{equation}
\label{eq:tube}
K_N(\mu)=\Phi\Big(\tfrac{\varepsilon_N-\Delta(\mu)}{s(\mu)}\Big)
-\Phi\Big(\tfrac{-\varepsilon_N-\Delta(\mu)}{s(\mu)}\Big)
\end{equation}
the Gaussian tube probability in the scalar case $r_0=1$, and its multivariate analogue
otherwise, we assume that on a fixed neighborhood $V$ of $E$ on which $\Omega$ is
nondegenerate,
\begin{equation}
\label{eq:lltscale}
\int_V\big|\Pr_\mu(E_N)-K_N(\mu)\big|\,\Pi(d\mu)=o(\varepsilon_N^{r_0}),
\end{equation}
and that $\Pr_\mu(E_N)\to0$ exponentially off $V$. Then $\Pi_{N,E}$ concentrates on $E$ at rate
$\varepsilon_N+N^{-1/2}$ and $\Pi_E=\Pi_\alpha$, the co-area conditional of
\eqref{eq:disint}, with no additional $\det\Omega(\mu)^{1/2}$ weighting beyond the
ordinary co-area disintegration.
\end{corollary}

\begin{proof}
Write $R_N(\mu)=\Pr_\mu(E_N)-K_N(\mu)$ on $V$, so that \eqref{eq:lltscale} reads
$\int_V|R_N|\,d\Pi=o(\varepsilon_N^{r_0})$. The remainder is controlled at the scale of the
normalization $\int\Pr_\mu(E_N)\Pi(d\mu)$, which is of order $\varepsilon_N^{r_0}$, or
$2\varepsilon_N=2N^{-1/2}$ in the scalar case with $\varepsilon_N=N^{-1/2}$. An additive
$o(1)$ would exceed it. Both $\varepsilon_N\to0$ and $s(\mu)\to0$, so for any $\omega_N\to\infty$ the weight is
negligible outside a neighborhood of $E$ of width
$\omega_N(\varepsilon_N+N^{-1/2})\to0$, which is the stated rate.

Write $\mu\leftrightarrow(a,v)$ for the disintegration \eqref{eq:disint},
$a=\langle g,\mu\rangle$ being the moment coordinate and $v$ a coordinate on the fiber. The key identity is that
for every fixed $s>0$, by Fubini,
\begin{equation}
\label{eq:coarea}
\int_{\Real}\Big\{\Phi\Big(\tfrac{\alpha+\varepsilon_N-a}{s}\Big)
-\Phi\Big(\tfrac{\alpha-\varepsilon_N-a}{s}\Big)\Big\}\,da=2\varepsilon_N,
\end{equation}
a density in the statistic integrating to one over the parameter. In the
multivariate case the same Fubini calculation on $\Real^{r_0}$ gives the tube volume
$(2\varepsilon_N)^{r_0}$, irrespective of correlations in $\Omega$. The right-hand
side does not depend on $s$.

$K_N$ depends on $\mu$ through $a$ and through
$s(\mu)=\{\Omega(\mu)/N\}^{1/2}$, and $\Omega$ is not a function of $a$ alone: it
varies along the fiber, running from $0.16$ to $0.96$ along $E$ in the example of
Section~\ref{sec:trap}. It cannot be frozen at a single value on $E$, but it can be
frozen in the normal direction, by continuity over a neighborhood of width
$\omega_N(\varepsilon_N+N^{-1/2})\to0$, giving $s(a,v)=s(\alpha,v)\{1+o(1)\}$
uniformly. For each fixed fiber coordinate $v$ the inner integral over $a$ is then
$2\varepsilon_N$ by \eqref{eq:coarea}, whatever value $s(\alpha,v)$ takes on that
fiber. The $\Omega$-dependence drops out fiber by fiber, and no covariance factor
survives.

Now test against a bounded continuous $f$. Integrating first in $a$ and then in $v$,
\eqref{eq:disint} and \eqref{eq:coarea} give
\begin{equation}
\label{eq:testfn}
\int_V f\,K_N\,d\Pi
=(2\varepsilon_N)^{r_0}h(\alpha)\int f\,d\Pi_\alpha+o(\varepsilon_N^{r_0}),
\end{equation}
where the error uses the assumed continuity of $h$, the weak continuity at $\alpha$
of the chosen version $a\mapsto\Pi_a$, and the continuity of $\Omega$ over that
neighborhood. Taking $f\equiv1$
in \eqref{eq:testfn} gives
$\int_V K_N\,d\Pi=\{h(\alpha)+o(1)\}(2\varepsilon_N)^{r_0}$, of exact order
$\varepsilon_N^{r_0}$ because $h(\alpha)>0$. Since
$\int_V|R_N|\,d\Pi=o(\varepsilon_N^{r_0})$, replacing $K_N$ by
$\Pr_\mu(E_N)=K_N+R_N$ perturbs both the numerator and denominator by only
$o(\varepsilon_N^{r_0})$. Consequently,
\[
\int f\,d\Pi_{N,E}
=\frac{\int f\,\Pr_\mu(E_N)\,\Pi(d\mu)}{\int \Pr_\mu(E_N)\,\Pi(d\mu)}
\longrightarrow\int f\,d\Pi_\alpha,
\]
the common leading factor $(2\varepsilon_N)^{r_0}h(\alpha)$ cancels. Since this
holds for every bounded continuous $f$, we obtain
$\Pi_E=\Pi_\alpha$, the co-area conditional.
\end{proof}

The zero-mass case deserves comment. When $\Pi(E)=0$, the notation
$\Pi(\cdot\mid E)$ has no intrinsic meaning until an approximating family and a
version on the null fiber are specified; this is the familiar Borel--Kolmogorov
issue. Here the weakly continuous co-area version fixes $\Pi_\alpha$, the entropy
tubes $\{\widehat{P}_N\in E_{\varepsilon_N}\}$ provide the approximating family, and
the theorem shows that the tubes select that version.

Large deviations alone do not identify this conditional. The logarithmic
approximation \eqref{eq:constraintLDP} omits the prefactor $a_N$, bounded in
\eqref{eq:typeprob} between $(N+1)^{-k}$ and $1$. On the relevant scale, this
prefactor carries the factor $\det\Omega^{-1/2}$ that cancels the tube Jacobian in
\eqref{eq:coarea}. Applying Laplace's method directly to
\eqref{eq:constraintLDP} misses the cancellation and instead gives
$\Pi_E(d\mu)\propto\det\Omega(\mu)^{1/2}\Pi(d\mu)$ on $E$. In
Section~\ref{sec:trap}, comparison with the exact finite-$N$ posterior yields the
co-area conditional.

Finally, condition \eqref{eq:accept} cannot be omitted if the prior is to remain
unreweighted along $E$. Let $\X=\{0,1,2\}$, set $E=\{Q:\Ex_QX=1\}$, and assign
equal prior mass to
\[
\mu_1=(0.25,\,0.50,\,0.25),\qquad \mu_2=(0.40,\,0.20,\,0.40).
\]
Both atoms lie in $E$, so $\Pi(E)=1$, but their variances are $0.5$ and $0.8$.
For $\varepsilon_N=N^{-1/2}$, their acceptance probabilities converge to
$2\Phi(1/\sigma_j)-1$, namely $0.843$ and $0.736$, and the limiting weights are
$0.534$ and $0.466$, rather than
$\tfrac12$ each: the tilt is inactive, as Proposition~\ref{prop:dichotomy}(ii) says, but
the prior is not preserved. A window with $\sqrt N\varepsilon_N\to\infty$ restores
\eqref{eq:accept} and with it the truncation.

\begin{remark}[Subsequential limits are free]
\label{rem:tight}
Hypothesis (i) of Theorem~\ref{thm:limit} is an assumption of weak convergence, but
half of it costs nothing. On a finite alphabet $\simplex(\X)$ is compact, so the
family $\{\Pi_{N,E}\}_{N\ge1}$ is automatically tight and every subsequence has a
further subsequence along which $\Pi_{N_k,E}\Rightarrow\Pi_E$ for some probability
measure $\Pi_E$; along that subsequence the conclusion of Theorem~\ref{thm:limit}
holds under (ii) and (iii) alone. This is the same compactness that supplies the Helly
limit in Theorem~\ref{thm:hs}, where the type lattice sits in $[0,1]$ and no tightness
question can arise either. What is not free is the identification of the limit point:
that Laplace's method applied to \eqref{eq:constraintLDP} selects the same
$I_E$-minimizing subfamily along every subsequence, which needs the $o_\mu(N)$ there
to be $o(N)$ uniformly on compacta. Hypothesis (i) may therefore be weakened to the
requirement that all subsequential limits induce the same pushforward
$(T_E)_\#\Pi_E$.
\end{remark}

The two regimes have different limiting effects. In the reachable case,
$I^\star=0$, the projection is asymptotically the identity. Sanov then contributes
through the finite-$N$ rate of Theorem~\ref{thm:master}, while any limiting
randomness comes from the constraint-selected prior; such randomness need not remain
when $\operatorname{supp}\Pi\cap E$ is a single point.

In the unreachable case, $I^\star>0$, the projection remains nontrivial, but the
constraint posterior concentrates on $\{I_E=I^\star\}$. Under the
affine, interiority, and rank conditions of
Proposition~\ref{prop:dichotomy}(iv), the limiting directing law is nondegenerate
if and only if $\operatorname{supp}\Pi_E$ meets more than one tilt orbit. Thus the
limit is both genuinely projected and nondegenerate precisely when $I^\star>0$ and
the limiting constraint posterior charges more than one orbit. The criterion
involves $\operatorname{supp}\Pi_E$, not the entire minimizing subfamily
$\{I_E=I^\star\}$: that subfamily may cross several orbits while the limiting
measure charges only one. None of the examples below has both properties at the
limit. At every finite $N$, however, componentwise conditioning and the
constraint-posterior update both enter the construction: Theorem~\ref{thm:master}
quantifies the former, while \eqref{eq:exactmix} displays the latter.

\begin{table}[H]
\centering
\footnotesize
\setlength{\tabcolsep}{3.5pt}
\caption{The dichotomy across the finite-alphabet examples.
$I^\star=\inf_{\operatorname{supp}\Pi}I_E$ decides whether the componentwise
projection remains nontrivial: in the reachable case it is asymptotically the
identity; in the unreachable case it remains nontrivial and, with the point-mass
priors of these two examples, the limiting mixing law is degenerate. In the
reachable rows, the corollaries rather than $I^\star$ identify the limit: the first
is a truncation (Corollary~\ref{cor:trunc}), and the second is a co-area conditional
(Corollary~\ref{cor:coarea}).}
\label{tab:dich}
\begin{tabular}{@{}l@{\ \ }c@{\ \ }l@{}}
\hline
Example & $I^\star$ & Limit \\
\hline
Beta$(2,2)$, $\mu\ge\alpha$ (\S\ref{sec:bernoulli})
 & $0$ & $\Pi(\cdot\mid E)$, $T_E=\mathrm{id}$; random \\
Dirichlet, mean $\alpha$ (\S\ref{sec:trap})
 & $0$ & co-area law, $T_E=\mathrm{id}$; random \\
$\Ber(0.2)$, mean $0.8$ (\S\ref{sec:bernoulli})
 & $>0$ & tilt acts; degenerate \\
Brandeis dice (\S\ref{sec:dice})
 & $>0$ & Jaynes tilt; degenerate \\
\hline
\end{tabular}
\end{table}

\subsection{Finite-horizon prediction}
\label{sec:finitepred}

Returning to the finite horizon, the constrained model has its own predictive
conditional. Let $J_N,I_N$ be disjoint with fixed cardinalities $t,m$. For
configurations $x_J,y_I$ whose limiting denominator is positive,
\begin{equation}
\label{eq:finiteratio}
\begin{aligned}
&\Qbb_N(X_{I_N}=y_I\mid X_{J_N}=x_J)\\
&\ \longrightarrow\
\frac{\int \prod_{i\in I_N}\Pstar_\mu(y_i)\prod_{j\in J_N}\Pstar_\mu(x_j)\,\Pi_E(d\mu)}
{\int \prod_{j\in J_N}\Pstar_\mu(x_j)\,\Pi_E(d\mu)},
\end{aligned}
\end{equation}
and conditionally on $\mu$,
$\Qbb_{N,\mu}(X_{I_N}=y_I\mid X_{J_N}=x_J)\to\prod_{i\in I_N}\Pstar_\mu(y_i)$. The
marginal statement follows by taking ratios of the $(t+m)$- and $t$-coordinate
limits in Theorem~\ref{thm:limit}; for each fixed $\mu$ satisfying the component
assumptions, the conditional-on-$\mu$ statement follows by the same ratio argument
applied to Theorem~\ref{thm:master}. We call \eqref{eq:finiteratio} the
\emph{constraint-conditioned predictive}, the predictive distribution within the
active constrained horizon. It is a valid conditional under $\Qbb_N$, and not an
ordinary predictive under the original unconstrained law.

\begin{remark}[Switching off the constraint]
\label{rem:switchoff}
There is a useful contrast with restarting the original model. Suppose the
$N$-horizon has been conditioned on $E_N$ and is then regarded as complete. If
one appends a new draw from the baseline component without imposing another
constraint, conditional independence gives
$\mathcal{L}(X_{N+1}\in A\mid E_N,\mu)=\mu(A)$. That procedure extends a completed
constrained sample by an unconstrained observation. It is not the
constrained-horizon experiment studied here, which enlarges the horizon and
imposes $E_{N+m}$ on the enlarged system.
\end{remark}

\subsection{Forward and inverse relative-entropy geometry}
\label{sec:invsanov}

Relative entropy appears at three different levels of the construction. To see
the first and the third directly, write the finite-alphabet likelihood as
\begin{equation}
\label{eq:likelihood}
\mu^{\otimes N}(x_{1:N})=\exp\{-NH(\widehat{P}_N)-N\KL{\widehat{P}_N}{\mu}\},
\end{equation}
the exact finite-$N$ form of the asymptotic equipartition property
\citep{shannon,mcmillan,breiman}. The entropy term, free of $\mu$, cancels in
Bayes' rule, so
\begin{equation}
\label{eq:invsanov}
\Pi(d\mu\mid x_{1:N})\ \propto\ e^{-N\KL{\widehat{P}_N}{\mu}}\,\Pi(d\mu).
\end{equation}

\begin{proposition}[Inverse-Sanov representation of the posterior]
\label{prop:invsanov}
On a finite alphabet, the de Finetti posterior over the directing measure is the
Gibbs measure \eqref{eq:invsanov}, at inverse temperature $N$ with energy
$\mu\mapsto\KL{\widehat{P}_N}{\mu}$.
\end{proposition}

Equation~\eqref{eq:invsanov} is classical. Its role here is to place three uses of
the same divergence side by side. Forward sampling assigns the type $Q$ the cost
$\KL{Q}{\mu}$. Conditioning selects $\Pstar_\mu$ through that forward cost.
Ordinary Bayes reverses the arguments and assigns the directing measure the cost
$\KL{\widehat{P}_N}{\mu}$. Under event conditioning, the effective energy is instead
$I_E(\mu)$. We use ``inverse Sanov'' only for this posterior Gibbs geometry. It
does not mean that a posterior large-deviation principle follows by formally
reversing every forward argument. Such a principle is established by
\citet{ganeshoconnell1999,ganeshoconnell}, while \citet{maccipiccioni} show why
the reversal can fail for exponential families.

\section{Thermodynamic reading and local fluctuations}
\label{sec:thermo}

The language of statistical mechanics is useful here, provided the reference
law $\mu$ is kept visible.

\emph{Microcanonical horizon.} Conditioning on
$\{\widehat{P}_N\in E_{\varepsilon_N}\}$ discards configurations whose empirical
moments fall outside the prescribed window. The remaining configurations retain
their weights under $\mu^{\otimes N}$; they are equiprobable only for a uniform
baseline. Thus the ensemble is microcanonical, or approximately so, relative to
$\mu$. The restriction belongs to the horizon as a whole, not to any one
coordinate.

\emph{Canonical local equilibrium.} The product $I$-projection
$(\Pstar_\mu)^{\otimes m}$ supplies the corresponding canonical law. Its Lagrange
multipliers tilt the baseline just enough to reproduce the prescribed moments.
Theorem~\ref{thm:master} makes the local equivalence precise: a fixed block from
the globally constrained horizon is close to independent draws from
$\Pstar_\mu$, with an explicit finite-$N$ error that does not depend on where the
block is located.

\emph{Free energy.} By the method of types
$\Pr_\mu(\widehat{P}_N=Q)=a_N(Q)\,e^{-N\KL{Q}{\mu}}$ with $(N+1)^{-k}\le a_N(Q)\le1$,
and
\begin{equation}
\label{eq:freeenergy}
\KL{Q}{\mu}=-H(Q)-\sum_x Q(x)\log\mu(x),
\end{equation}
The summation term in \eqref{eq:freeenergy} is the energy relative to the
baseline, while $H(Q)$ is the entropy. Their difference $\KL{Q}{\mu}$ is the
intensive free-energy cost of the macrostate $Q$.
Minimizing that cost over the admissible set gives the $I$-projection
$\Pstar_\mu$.

\emph{Entropy-energy balance.} The bad-type mass in the proof obeys
$\beta_N(\delta)\le C(\mu,E)\exp\{2k\log(N+1)-2N\delta^2\}$. The logarithmic term
accounts for the at most $(N+1)^k$ types and the corresponding normalization of the
conditioning denominator; the extensive term $2N\delta^2$ is the
relative-entropy barrier.
Balancing them gives $N\delta_N^2\asymp k\log N$, hence
$\delta_N\asymp\sqrt{(k\log N)/N}$, the tuning of
Corollary~\ref{cor:tuning}. This split radius, distinct from the window radius
$\varepsilon_N$, is where the combinatorial entropy of the type lattice is balanced
by the divergence cost of leaving the projection.

\emph{Relation to PAC-Bayes.} The same architecture governs Gibbs-posterior and
PAC-Bayesian concentration \citep{catoni}. PAC-Bayes bounds the excess risk of a
randomized predictor at a square-root complexity-to-sample-size scale, with
complexity measured by relative entropy between posterior and prior.
Corollary~\ref{cor:tuning} has the analogous scale
$\sqrt{(k\log N)/N}$: the logarithmic complexity of the type lattice is divided by
the horizon, while excess relative entropy from the constrained equilibrium
provides the barrier. Moreover, both the Gibbs posterior of PAC-Bayes and the
conditional type law $w_N$ in \eqref{eq:splitset} are exponential tilts of a
reference measure by an empirical cost. The correspondence is structural, not a
claim of a generalization bound: here the state space consists of empirical types,
the cost is excess relative entropy from canonical equilibrium, and the conclusion
is a local equivalence of ensembles.

\subsection{Local fluctuations on the
\texorpdfstring{$N^{-1/2}$}{inverse-square-root N} scale}
\label{sec:fluctuations}

The master inequality is a global concentration statement, controlling the
probability of leaving a fixed or slowly shrinking neighborhood of $\Pstar$. Inside
that neighborhood, on the $N^{-1/2}$ scale, there is a complementary local
refinement, which we record for the regular affine case; the global
Gibbs-conditioning theorem does not by itself supply it.

\begin{remark}[Local Gaussian fluctuations; companion paper]
\label{rem:fluct}
Under exact feasibility on the type lattice and an affine local slice with
no active inequality, the conditional law of $\sqrt N(\widehat{P}_N-\Pstar)$ is
asymptotically Gaussian on the constraint tangent space, with precision the reduced
Hessian of $\KL{\cdot}{\mu}$ at $\Pstar$. Consequently
$2N\,\KL{\widehat{P}_N}{\Pstar}$ is asymptotically chi-square with $k-1-r_0$ degrees of
freedom. Here $r_0$ is the number of equality constraints; no inequality is active
on the slice, whereas $r$ counts constraints of both kinds.
\end{remark}

Neither the master inequality nor a shrinking window supplies this refinement. Any
admissible window satisfies $N\varepsilon_N>c_{\mathrm{rnd}}$ and is therefore thicker
than the $O(1/N)$ mesh of the type lattice, so the feasible types acquire a normal
component and a two-scale local limit theorem would be required, resolving $N^{-1/2}$
tangential fluctuations together with a normal scale determined jointly by the shrinking
window and the $O(1/N)$ lattice mesh. Exact feasibility collapses the feasible set
onto a single tangent slice, which is what makes the single-scale Gaussian available.
The companion paper \citep{polsonzantedeschi} develops the exact-feasibility
Gaussian theorem, its chart geometry, and the associated estimation theory, while
leaving the thickened case open. Thus Theorem~\ref{thm:master} and
Remark~\ref{rem:fluct} are complementary: the former gives global concentration
under a window, the latter local geometry without one.

Two distinct local Gaussian phenomena appear, in opposite arguments of the
divergence. Forward conditional Sanov gives Gaussian fluctuations of the empirical
type around the $I$-projection on the constraint tangent space
(Remark~\ref{rem:fluct}). Inverse Sanov gives Gaussian posterior fluctuations of
the directing law around the empirical distribution: for a true interior component
$\mu_0$ and a prior density positive and continuous at $\mu_0$, the ordinary
finite-dimensional Bernstein--von Mises theorem \citep[Theorem~10.1]{vandervaart} gives, in any
smooth $(k-1)$-dimensional chart,
\[
\mathcal{L}\big(\sqrt N(\mu-\widehat{P}_N)\mid X_{1:N}\big)\ \Longrightarrow\
N_{k-1}\big(0,\Fish(\mu_0)^{-1}\big),
\]
the ambient covariance being
$\operatorname{diag}(\mu_0)-\mu_0\mu_0^\top$ restricted to the simplex tangent
space. Both arise from the local quadratic expansion of relative entropy, in
opposite arguments: the empirical type fluctuates around the projected law, and the
posterior directing law fluctuates around the empirical type.

\section{Examples}
\label{sec:examples}

\subsection{Bernoulli components}
\label{sec:bernoulli}

Example~\ref{ex:bernoulli} is the simplest case: with $\mu=\Ber(0.2)$ and the horizon
mean constrained near $0.8$, every fixed target inside the horizon is
asymptotically $\Ber(0.8)$, the $I$-projection, regardless of its temporal label.

The binary case is small enough to evaluate exactly. Conditionally on a type the
block is hypergeometric, so \eqref{eq:hs} gives the constrained block law in closed
form and its distance to $(\Pstar)^{\otimes m}$ needs no simulation.
Table~\ref{tab:bern} separates the two contributions of
Theorem~\ref{thm:master} for $m=5$ and $\varepsilon_N=N^{-1/2}$: the split-set term,
which collects the type concentration and the finite-population coupling, and the
window-transfer term $\TV{(P^\bullet)^{\otimes m}-(\Pstar)^{\otimes m}}$.

\begin{table}[t]
\centering
\footnotesize
\caption{Exact total variation for the constrained block law. $\mu=\Ber(0.2)$,
horizon mean $0.8$, $m=5$, $\varepsilon_N=N^{-1/2}$. All three entries are exact
distances: split-set is $\TV{L_{N,I_N,\mu}-(P^\bullet)^{\otimes m}}$, window is
$\TV{(P^\bullet)^{\otimes m}-(\Pstar)^{\otimes m}}$, total is
$\TV{L_{N,I_N,\mu}-(\Pstar)^{\otimes m}}$. They are related by the triangle
inequality, not by addition; the window term is the binding contribution.}
\label{tab:bern}
\begin{tabular}{@{}r@{\ \ }c@{\ \ }c@{\ \ }c@{\ \ }c@{}}
\hline
$N$ & $\varepsilon_N$ & split-set & window & total \\
\hline
$25$ & $0.200$ & $0.045$ & $0.400$ & $0.410$ \\
$100$ & $0.100$ & $0.013$ & $0.209$ & $0.209$ \\
$400$ & $0.050$ & $0.003$ & $0.104$ & $0.104$ \\
$800$ & $0.035$ & $0.001$ & $0.074$ & $0.073$ \\
$1600$ & $0.025$ & $0.001$ & $0.052$ & $0.052$ \\
\hline
\end{tabular}
\end{table}

At every horizon the split-set term is smaller than its coupling budget
$m(m-1)/(2N)$ by a factor between six and twelve: it is $0.045$ against $0.400$
at $N=25$, and $0.003$ against $0.025$ at $N=400$. The window width is therefore
the binding contribution.

The equality constraint here is affine with $r_0=1$ on an alphabet of size $k=2$, so
$k-1-r_0=0$ and Proposition~\ref{prop:dichotomy}(iv) forces the projected directing
measure to be degenerate: $E$ is the singleton $\{\Ber(0.8)\}$ and $\Pstar_\mu=\Ber(0.8)$
for every $\mu$. Numerically, at $N=800$ with $m=2$, the constrained block law agrees to within
$0.005$ in total variation across the baselines $\Ber(0.05)$, $\Ber(0.2)$, $\Ber(0.5)$,
and $\Ber(0.7)$, each sitting within $0.056$ of $\Ber(0.8)^{\otimes m}$: no prior
survives.

The inequality constraint $E=\{Q:Q(1)\ge\alpha\}$ is a different branch, and clause
(iv) does not apply to it: $E$ is not affine, and $T_E$ is the identity on
$\{\mu\ge\alpha\}$ and constant on $\{\mu<\alpha\}$, so it is neither an orbit quotient
nor injective. It is instead the reachable branch, and
Corollary~\ref{cor:trunc} applies to it directly: $\Pi(\partial E)=\Pi(\{\alpha\})=0$
in $\simplex(\X)$, so the acceptance condition \eqref{eq:accept} holds. Because the
$\mathrm{Beta}(2,2)$ prior has topological support $[0,1]$, this conclusion comes
from the corollary rather than Proposition~\ref{prop:dichotomy}. With
$\Pi=\mathrm{Beta}(2,2)$ and $\alpha=0.6$ the set $E$ carries $35.2\%$ of the prior, so
$\Pi_E=\Pi(\cdot\mid\mu\ge\alpha)$, with $\Ex[\mu\mid\mu\ge\alpha]=0.745$ and standard
deviation $0.097$; the finite-horizon posterior means rise to it, from $0.619$ at $N=25$
to $0.735$ at $N=3200$. The directing measure is genuinely random. On
$\{\mu\ge\alpha\}$, however, the projection is the identity, so the limiting randomness
is the truncated prior, while the componentwise tilt contributes only at finite
$N$ through the bound in Table~\ref{tab:bern}.

\subsection{A null constraint set: the co-area conditional}
\label{sec:trap}

For a regular prior and a prior-null equality constraint,
Corollary~\ref{cor:coarea} gives the relevant limit; in this setting the
large-deviation form \eqref{eq:constraintLDP} does not by itself identify the
conditional.
Specifically, take $\X=\{1,2,3\}$, $g(x)=x$, the single equality constraint
$E=\{Q:\Ex_Q[X]=\alpha\}$ with $\alpha=2.2$, and the flat prior
$\Pi=\mathrm{Dirichlet}(1,1,1)$. Then $\Pi(E)=0$, because $E$ is a segment in the
simplex, while $I_E$ vanishes on $E$. The hypotheses of
Corollary~\ref{cor:coarea} hold: the pushforward of $\Pi$ under
$\mu\mapsto\Ex_\mu[X]$ has a continuous positive density at $\alpha=2.2$.

The segment is $\mu(u)=(u-0.2,\,1.2-2u,\,u)$ for $u\in[0.2,0.6]$, with
$\sigma^2(\mu(u))=\operatorname{Var}_\mu(X)=2u-0.24$, running from $0.16$ at the die
$(0,0.8,0.2)$ to $0.96$ at the die $(0.4,0,0.6)$. The constraint map is affine and the
prior is flat, so the co-area disintegration is uniform in $u$ and predicts
$\Ex[\mu_3\mid E]=0.4$ exactly. Laplace's method applied to \eqref{eq:constraintLDP}
alone predicts instead the $\sigma$-reweighted law $\Pi_E\propto\sigma(\mu)\,\Pi$, giving
$0.4255$. The two differ by $2.6$ percentage points in a predictive probability, since
$\mu_3$ is the limiting $\Pr(X_{t+1}=3)$.

The finite-$N$ constraint posterior settles it. The type law is trinomial, so
$\Pr_\mu(E_N)$ is a finite sum and $\Pi_{N,E}$ can be evaluated exactly, without
simulation, under either the window $E_{\varepsilon_N}$ or exact conditioning on
$\{\bar X_N=2.2\}$. Table~\ref{tab:trap} reports it.

\begin{table}[H]
\centering
\footnotesize
\caption{$\Ex[\mu_3\mid E_N]$ under the exact constraint posterior.
$\X=\{1,2,3\}$, $\Ex_Q[X]=2.2$, $\Pi=\mathrm{Dirichlet}(1,1,1)$. Co-area predicts
$0.400$; the $\sigma$-reweighted law predicts $0.4255$.}
\label{tab:trap}
\begin{tabular}{@{}r@{\quad}c@{\quad}c@{}}
\hline
$N$ & window, $\varepsilon_N=N^{-1/2}$ & exact, $\{\bar X_N=2.2\}$ \\
\hline
$100$ & $0.394$ & $0.398$ \\
$400$ & $0.398$ & $0.400$ \\
$1600$ & $0.400$ & $0.400$ \\
\hline
\end{tabular}
\end{table}

Both schemes converge to $0.400$ from below at rate $O(1/N)$, and the
$\sigma$-reweighting is excluded. The reason is the cancellation in the proof of
Corollary~\ref{cor:coarea}: the polynomial prefactor omitted by
\eqref{eq:constraintLDP} carries $\det\Omega^{-1/2}$, which cancels the
Jacobian $\det\Omega^{1/2}$ contributed by the Gaussian tube. A large-deviation rate
identifies where the posterior concentrates; it does not by itself identify the
conditional it concentrates to.

\subsection{The Brandeis dice problem}
\label{sec:dice}

For each active horizon $N$, impose the empirical-average constraint
$\sum_{i=1}^{6} i\,\PN(i)\approx 4.5$, against the fair value $3.5$. The baseline is
uniform, so the $I$-projection is genuinely a maximum-entropy tilt,
\[
p_i=\frac{e^{\lambda i}}{\sum_{j=1}^6 e^{\lambda j}},\qquad \lambda\approx0.37105,
\]
giving $p\approx(0.054,0.079,0.114,0.165,0.240,0.347)$ with
$H\approx1.614<\log6$. See \citet{jaynes} for the maximum-entropy analysis and
\citet{seidenfeld} for a critical discussion of its Bayesian predictive
interpretation. Any fixed finite collection of
throws in that horizon has a local law converging to the product of this tilt,
including throws labeled future relative to any finite cut; when the baseline die
is fixed and uniform, the limiting constrained process is i.i.d.\ from the Jaynes
tilt. The uniform die has mean $3.5\neq4.5$, so $I^\star>0$: this is the unreachable
branch of Proposition~\ref{prop:dichotomy}(iii), with a point-mass prior, hence
degenerate (Table~\ref{tab:dich}). Under uncertain die laws, $\PiE=(\mu\mapsto\Pstar_\mu)_\#\Pi_E$ is the induced
prior over constrained die laws and sequential prediction is the Bayesian mixture.

\subsection{Gaussian scale mixtures}
\label{sec:gaussian}

Under the usual exponential-integrability and regularity conditions, the continuous
Gibbs conditioning principle extends the same structure beyond finite alphabets
\citep{leonardnajim}. We use Gaussian scale mixtures as an illustration.

The Gaussian components arise from symmetry rather than assumption.
\citet{kingman} shows that if every finite sample from an infinite sequence has a
spherically symmetric joint law, then conditionally on a single random variable $V$
the coordinates are i.i.d.\ $N(0,V)$. The argument is a characteristic-function one:
spherical symmetry forces the characteristic function of a finite block to depend on
$t$ only through $\lVert t\rVert$; Schoenberg's theorem \citep{schoenberg}
identifies the admissible radial forms as Gaussian scale mixtures, with mixing
variable $V$.
\citet{smith} gives the centered variant, with a random mean as well, and
\citet{ressel} gives an analytic account. The result is the natural continuous
counterpart of the scalar-type construction in Section~\ref{sec:hs}: compactness of
$[0,1]$ and Helly selection produce the mixing law there, whereas spherical
symmetry and Schoenberg's theorem produce it here. In particular, the theorem
supplies the prior over Gaussian directing components on which the constraint acts.

Although Proposition~\ref{prop:dichotomy} is stated for a finite alphabet, the same
interpretation applies. Suppose the prior over $(M,V)$
charges every neighborhood of a component with the prescribed moments, so that
$E\cap\operatorname{supp}\Pi\neq\varnothing$ and $I^\star=0$. Reachability is
relative to $\operatorname{supp}\Pi$; the mere existence of a Gaussian with the
prescribed moments is not enough. Then $T_E$ is the identity on $E$, and the
surviving heterogeneity depends on how the constraint reweights the prior. Under
mean-only conditioning, the $I$-projection of $N(M,V)$ is $N(\theta,V)$, which
retains the variance. The limit therefore mixes $N(\theta,V)$ over a conditional
law of $V$ and is nondegenerate exactly when that conditional law is.

The form of the conditional law depends on the regime. If $\Pi$ has a regular
density in $(M,V)$, the
constraint is $\Pi$-null and the same normal-direction cancellation as in
Corollary~\ref{cor:coarea} applies: for fixed $V$,
\begin{equation}
\label{eq:gaussfubini}
\int_{\Real}\Pr_{M,V}\big(|\bar X_N-\theta|\le\varepsilon_N\big)\,dM=2\varepsilon_N
\end{equation}
for every $V$, so the limit is the co-area conditional at $M=\theta$ carrying \emph{no}
additional variance factor. A larger $V$ does lower acceptance at the exact component
$M=\theta$, but it widens the range of nearby $M$ that can produce the event, and the
two effects cancel exactly. If instead $\Pi$ places positive mass on $\{M=\theta\}$, the
constraint is not $\Pi$-null, \eqref{eq:gaussfubini} does not apply, and a window of
order $N^{-1/2}$ reweights those atoms through their componentwise acceptance
probabilities, as in Section~\ref{sec:dichotomy}. This is a continuous analogue of
Corollary~\ref{cor:coarea}, not an application of it, because the corollary is stated
for finite alphabets. Under mean-and-variance conditioning, the projection is
$N(\theta,v)$ for every component; both parameters are pinned, and the mixture
collapses to the single product law $N(\theta,v)^{\otimes m}$. A finite-block
continuous theorem with rates in the present mixture-level, shrinking-window
generality remains open.

\section{Discussion}
\label{sec:discussion}

We treat the empirical-moment restriction as a primitive global constraint imposed
separately on each active horizon. The target is local while the constraint is
global; exchangeability makes the target's location irrelevant, and Gibbs
conditioning converts the constrained component law into a local product
$I$-projection. After mixing and passing to the local limit, the projected
components define the law of a random directing measure for a genuine exchangeable
process, under which prediction is
ordinary Bayesian prediction: minimum relative entropy within components, Bayesian
mixing across projected components. Proposition~\ref{prop:dichotomy} distinguishes
the two limiting regimes. Under a reachable constraint the projection becomes
asymptotically the identity; under an unreachable constraint the constraint posterior
concentrates on the $I_E$-minimizing subfamily. The conditional selected within that
subfamily is determined by Corollaries~\ref{cor:trunc} and~\ref{cor:coarea} under
their respective hypotheses. Componentwise conditioning and the mixing update both
enter at every finite horizon.

The present paper and its companion address opposite sides of the mixture. Here the
question is which directing law results from a primitive constraint. The companion
instead studies the rate and geometry with which a block approaches a single
$I$-projection, together with the inference supported by that localization.

Several limitations remain. The convergence theorem assumes a finite alphabet and
(A1)--(A5). For the broader exchangeability framework on general state spaces, see
\citet{kallenberg}; a finite-block continuous analogue with rates in the present
mixture-level, shrinking-window generality remains open.
The rate obtained here is an upper bound, with sharper
regular-case results available from \citet{diaconisfreedman1988}. Finally,
Proposition~\ref{prop:dichotomy} locates the support of $\Pi_E$ conditional on the
existence of a weak limit. Remark~\ref{rem:tight} makes subsequential existence
automatic on a finite alphabet, but uniqueness remains to be established. Moreover,
as Corollary~\ref{cor:coarea} and Section~\ref{sec:trap} show, a pointwise
large-deviation rate does not identify the conditional within the minimizing
subfamily. That identification requires local central-limit information, and the
resulting weights depend on the window.

\begin{appendix}

\section{Notation}
\label{app:notation}

Table~\ref{tab:notation} collects the recurring symbols, each with the equation or
result that defines it.

\begin{table}[H]
\centering
\footnotesize
\setlength{\tabcolsep}{4pt}
\caption{Notation, and where each symbol is defined.}
\label{tab:notation}
\begin{tabular}{@{}l@{\ \ }l@{\ \ }l@{}}
\hline
Symbol & Meaning & Defined \\
\hline
$\X,\ k$ & alphabet and its size & \S\ref{sec:types}\\
$N,\ m$ & horizon, block size & \S\ref{sec:local}\\
$I_N$ & target index set, $|I_N|=m$ & \eqref{eq:locfree}\\
$\PN$ & empirical measure of $X_{1:N}$ & \eqref{eq:eventN}\\
$\TN$ & types with denominator $N$ & \S\ref{sec:types}\\
$\mu,\ \Pi$ & directing measure, prior & Thm~\ref{thm:definetti}\\
$F_N$ & law of the scalar type, $k=2$ & \S\ref{sec:hs}\\
$E$ & moment-constraint set & \eqref{eq:momentset}\\
$E_{\varepsilon_N}$ & moment window, radius $\varepsilon_N$ & \eqref{eq:window}\\
$E_N$ & horizon event $\{\PN\in E_{\varepsilon_N}\}$ & \eqref{eq:eventN}\\
$\Qbb_N$ & constrained-horizon law & \eqref{eq:QN}\\
$\Qbb_{N,\mu}$ & constrained component law & \eqref{eq:locfree}\\
$\Pi_{N,E}$ & constraint posterior & \eqref{eq:constraintpost}\\
$\Pi_E$ & weak limit of $\Pi_{N,E}$ & Thm~\ref{thm:limit}\\
$\KL{Q}{\mu}$ & relative entropy & \S\ref{sec:background}\\
$\TV{\cdot}$ & total variation & \S\ref{sec:background}\\
$\Pstar_\mu$ & $I$-projection of $\mu$ on $E$ & \eqref{eq:iproj}\\
$P^\bullet$ & $I$-projection on the window & \S\ref{sec:masterthm}\\
$I_E(\mu)$ & $\inf_{Q\in E}\KL{Q}{\mu}$ & \eqref{eq:constraintLDP}\\
$I^\star$ & $\inf_{\operatorname{supp}\Pi}I_E$; reachability & Prop~\ref{prop:dichotomy}\\
$\Omega(\mu)$ & $\operatorname{Cov}_\mu(g)$, moment covariance & \eqref{eq:tube}\\
$T_E$ & the map $\mu\mapsto\Pstar_\mu$ & Thm~\ref{thm:limit}\\
$\PiE$ & $(T_E)_\#\Pi_E$, directing law & Thm~\ref{thm:limit}\\
$Q^\star$ & random directing measure & Cor~\ref{cor:bayespred}\\
$w_N(Q)$ & conditional type weight & \S\ref{sec:masterthm}\\
$G_{N,\delta},B_{N,\delta}$ & good and bad types & \eqref{eq:splitset}\\
$\beta_N$ & bad-set mass & Prop~\ref{prop:split}\\
$\delta_N$ & split radius & Cor~\ref{cor:tuning}\\
$c_{\mathrm{rnd}}$ & rounding constant & \S\ref{sec:windows}\\
\hline
\end{tabular}
\end{table}

\section{Proofs}
\label{app:proofs}

Throughout, $P>0$ satisfies (A1)--(A5), $F_N:=E_{\varepsilon_N}$ and
$P^\bullet:=\Pstar_{\varepsilon_N}$.

\begin{lemma}[Window stability]
\label{lem:stab}
Under \textup{(A1)--(A5)}, $\TV{P^\bullet-\Pstar}\le C_0\varepsilon_N$ and
$\TV{\Pstar_{\varepsilon_N-c_{\mathrm{rnd}}/N}-P^\bullet}\le C_0 c_{\mathrm{rnd}}/N$.
\end{lemma}

\begin{proof}
Both are instances of the pairwise bound (A5) with $\varepsilon'=0$ and
$\varepsilon'=\varepsilon_N-c_{\mathrm{rnd}}/N$ respectively. The sufficient
conditions are the constraint qualifications of Section~\ref{sec:assumptions}: the
window relaxes each active bound by $O(\varepsilon)$, and Robinson strong
regularity of the KKT system makes the primal solution Lipschitz in the
right-hand sides, with the second-order condition supplied by positive
definiteness of the Hessian of $\KL{\cdot}{P}$ at the interior solution.
\end{proof}

\begin{lemma}[Inner-window type and denominator]
\label{lem:denom}
Under \textup{(A1)--(A5)} with $N\varepsilon_N>c_{\mathrm{rnd}}$, there is a type
$Q_N\in F_N\cap\TN$ with $\KL{Q_N}{P}\le\KL{P^\bullet}{P}+O(1/N)$; consequently, for
all sufficiently large $N$,
\[
\begin{aligned}
\sum_{Q'\in F_N\cap\TN}\Pr_P(\widehat{P}_N=Q')
&\ge (N+1)^{-k}\\[-2pt]
&\quad{}\times e^{-N\KL{P^\bullet}{P}}\,e^{-O(1)} .
\end{aligned}
\]
\end{lemma}

\begin{proof}
Let $P^{\bullet\bullet}=\Pstar_{\varepsilon_N-c_{\mathrm{rnd}}/N}$; by
Lemma~\ref{lem:stab}, $\TV{P^{\bullet\bullet}-P^\bullet}=O(1/N)$. Round
$P^{\bullet\bullet}$ to a type $Q_N$ with $\TV{Q_N-P^{\bullet\bullet}}=O(1/N)$, so
each moment moves by at most $c_{\mathrm{rnd}}/N$ and the inner-window margin gives
$Q_N\in F_N$. Interiority (A4) and $P>0$ (A2) make $\KL{\cdot}{P}$ Lipschitz near
$Q_N$ for large $N$, whence
$\KL{Q_N}{P}\le\KL{P^{\bullet\bullet}}{P}+O(1/N)\le\KL{P^\bullet}{P}+O(1/N)$. The
lower estimate in \eqref{eq:typeprob} at $Q_N$ gives the bound, with
$e^{-O(1)}=C(P,E)^{-1}$.
\end{proof}

\begin{lemma}[Divergence-gap concentration]
\label{lem:gap}
Under \textup{(A1)--(A5)}, for $\delta>0$,
\[
\begin{aligned}
&\sum_{\substack{Q\in F_N\cap\TN\\ \TV{Q-P^\bullet}>\delta}}
\Pr_P(\widehat{P}_N=Q\mid \widehat{P}_N\in F_N)\\
&\qquad\le C(P,E)(N+1)^{2k}e^{-2N\delta^2},
\end{aligned}
\]
with $C(P,E)$ from Lemma~\ref{lem:denom}.
\end{lemma}

\begin{proof}
By \eqref{eq:typeprob} the numerator of each conditional weight is at most
$e^{-N\KL{Q}{P}}$, and by Lemma~\ref{lem:denom} the denominator is at least
$C(P,E)^{-1}(N+1)^{-k}e^{-N\KL{P^\bullet}{P}}$, so
\[
w_N(Q)\le C(P,E)(N+1)^k\,e^{-N(\KL{Q}{P}-\KL{P^\bullet}{P})}.
\]
Since $F_N$ is convex with $I$-projection $P^\bullet$, \eqref{eq:klgap} gives
$\KL{Q}{P}-\KL{P^\bullet}{P}\ge 2\TV{Q-P^\bullet}^2$, which exceeds $2\delta^2$ on
the sum; there are at most $(N+1)^k$ such types.
\end{proof}

\begin{proof}[Proof of Theorem~\ref{thm:master}]
By exchangeability of the constrained component law \eqref{eq:locfree} it suffices
to take $I_N=\{1,\dots,m\}$. Decompose using \eqref{eq:hs} as
\[
L_{N,I_N,P}=\sum_{Q\in F_N\cap\TN}w_N(Q)\,\mathcal{L}_Q,
\]
with $\mathcal{L}_Q=\mathcal{L}_P(X_{1:m}\mid \widehat{P}_N=Q)$ and
$w_N(Q)=\Pr_P(\widehat{P}_N=Q\mid \widehat{P}_N\in F_N)$, and split the types at
$\TV{Q-P^\bullet}\le\delta$.
For good types, \eqref{eq:coupling} and
$\TV{Q^{\otimes m}-(P^\bullet)^{\otimes m}}\le m\TV{Q-P^\bullet}\le m\delta$ give
$\TV{\mathcal{L}_Q-(P^\bullet)^{\otimes m}}\le m(m-1)/(2N)+m\delta$; averaging
preserves this. By Lemma~\ref{lem:gap} the bad mass
$\Pr_P(\widehat{P}_N\in B_{N,\delta}\mid E_N)\le\beta_N:=C(P,E)(N+1)^{2k}e^{-2N\delta^2}$, and
since each $\TV{\cdot}\le1$ its contribution is at most $\beta_N$; this is the
split-set bound of Proposition~\ref{prop:split}. Finally
$\TV{(P^\bullet)^{\otimes m}-(\Pstar)^{\otimes m}}\le m\TV{P^\bullet-\Pstar}\le
C_0 m\varepsilon_N$ by Lemma~\ref{lem:stab}, and the triangle inequality gives
\eqref{eq:master}. For Corollary~\ref{cor:tuning}, $\delta_N$ gives
$2N\delta_N^2=2(k+1)\log(N+1)$, so $\beta_N\le C(P,E)(N+1)^{-2}$.
\end{proof}

\end{appendix}

\end{document}